\documentclass[12pt]{article}
\usepackage[margin=1.25in]{geometry}
\usepackage{amssymb}
\usepackage{amsmath}
\usepackage{amsthm}
\usepackage{hyperref}
\usepackage{mathrsfs}
\usepackage{textcomp}
\hypersetup{hidelinks, pdftitle={Stability of separable solutions to a fast diffusion system}, pdfauthor={Xuzhou Yang}}

\usepackage{verbatim}
\usepackage{esint}

\usepackage{tikz}

\makeatletter

\newdimen\bibspace
\renewenvironment{thebibliography}[1]{%
 \section*{\refname
       \@mkboth{\MakeUppercase\refname}{\MakeUppercase\refname}}%
     \list{\@biblabel{\@arabic\c@enumiv}}%
          {\settowidth\labelwidth{\@biblabel{#1}}%
           \leftmargin\labelwidth
           \advance\leftmargin\labelsep
           \itemsep\bibspace
           \parsep\z@skip     %
           \@openbib@code
           \usecounter{enumiv}%
           \let\p@enumiv\@empty
           \renewcommand\theenumiv{\@arabic\c@enumiv}}%
     \sloppy\clubpenalty4000\widowpenalty4000%
     \sfcode`\.\@m}
    {\def\@noitemerr
      {\@latex@warning{Empty `thebibliography' environment}}%
     \endlist}

\makeatother

\newtheorem{thm}{Theorem}[section]
\newtheorem{lem}[thm]{Lemma}
\newtheorem{prop}[thm]{Proposition}
\newtheorem{defn}[thm]{Definition}

\newtheorem{cor}[thm]{Corollary}
\newtheorem{rem}[thm]{Remark}

\def\XXint#1#2#3{{\setbox0=\hbox{$#1{#2#3}{\int}$}
  \vcenter{\hbox{$#2#3$}}\kern-.5\wd0}}

\newcommand{\bpf}{\begin{proof}}      \newcommand{\epf}{\end{proof}}
\newcommand{\bpr}{\begin{prop}}      \newcommand{\epr}{\end{prop}}
\newcommand{\brem}{\begin{rem}}      \newcommand{\erem}{\end{rem}}
\newcommand{\bthm}{\begin{thm}}      \newcommand{\ethm}{\end{thm}}
\newcommand{\blem}{\begin{lem}}      \newcommand{\elem}{\end{lem}}

\newcommand{\om}{\Omega}                \newcommand{\pa}{\partial}
           \newcommand{\ud}{\mathrm{d}}
\newcommand{\be}{\begin{equation}}      \newcommand{\ee}{\end{equation}}

\newcommand{\R}{\mathbb{R}}

\newcommand{\dist}{\mbox{dist}}
\newcommand{\normmm}[1]{{\left\vert\kern-0.25ex\left\vert\kern-0.25ex\left\vert #1
   \right\vert\kern-0.25ex\right\vert\kern-0.25ex\right\vert}}

\begin{document}

\title{\textbf{Stability of separable solutions to a fast diffusion system}
\bigskip}

\author{\large Xuzhou Yang\thanks{School of Mathematics, Yunnan Normal University, Kunming, China. Email: \texttt{2768620240@qq.com}.}}

\maketitle
\rmfamily
\begin{abstract}
We study a Dirichlet fast diffusion system with symmetric linear cross-loss terms on a bounded smooth domain. A coupled Rayleigh quotient captures the variational structure of the coupling and is nonincreasing along sufficiently regular positive solution pairs. Under the common extinction and regularity hypotheses, this monotonicity and the energy identity yield upper and lower bounds
\[
C^{-1}(T^*-t)^{1/(p-1)}
\le \|(u(t),v(t))\|_{L^{p+1}(\Omega)}
\le C(T^*-t)^{1/(p-1)}.
\]
In the Sobolev subcritical range, a Lyapunov functional for the rescaled system yields a sequence of times $t_j\uparrow T^*$ along which the rescaled pair converges strongly in $H_0^1(\Omega)\times H_0^1(\Omega)$ to a nonnegative classical stationary solution pair. The integral lower bound excludes the trivial limiting pair, while the coupling with $\alpha>0$ ensures that neither limiting component is identically zero. As a further property of the stationary system, we show that proportional components of a nontrivial nonnegative solution pair must coincide. We also construct a positive diagonal stationary solution pair by constrained minimization. Scalar upper comparison and a B\'enilan--Crandall-type estimate for smooth positive diagonal solutions provide complementary links with the scalar fast diffusion equation.
\end{abstract}


\smallskip


\bigskip

\tableofcontents

\section{Introduction}\label{sec:intro}
Let $ \Omega \subset \R^n$, $n \geq 1 $,  be a bounded smooth domain and $1<p< +\infty$ be a real number. We consider the nonlinear diffusion system under the Dirichlet boundary condition
\be\label{eq:fast diffusion system}
  \left\{     \begin{aligned}     \pa_t u^p & =\Delta  u-\alpha v \quad              \mbox{in }   \om \times (0,\infty),\\
 \pa_t v^p  &    =\Delta v -\alpha u\quad                  \mbox{in }   \om \times (0,\infty),\\
 u &=v=0  \quad  \mbox{on } \pa \om \times (0,\infty) ,\\
 \left( u(x,0),v(x,0) \right) &    = \left(u_0(x),v_0(x)\right)   \quad  \mbox{on }  \om \times \{ t=0 \},
\end{aligned}
\right.
\ee
Here $\alpha>0$ is fixed, $\Delta$ acts on the spatial variable, and the initial data $u_0,v_0\in C_0^1(\Omega)$ are assumed to be nonnegative. In the model introduced here, each component contributes a linear loss to the other. Coupling provides a way to describe interaction, as in the classical Lotka--Volterra equations \cite{lotka,volterra,2022CD-DLVS}. We refer to \eqref{eq:fast diffusion system} as a fast diffusion system because both diffusion terms have the power-law structure of the scalar Dirichlet fast diffusion equation:
\begin{equation}\label{eq:fde}
 \left\{
\begin{aligned}
& \partial_{t}u^{q}  =\Delta u \quad                           \mbox{in }   \Omega \times (0,\infty),\\
&u              =  0  \quad   \mbox{on } \partial \Omega \times [0,\infty).
\end{aligned}
\right.
\end{equation}
In the case of $1<q<+\infty$, \eqref{eq:fde} is called the fast diffusion equation. It arises in models of gas kinetics, plasma diffusion and thin liquid films driven by van der Waals forces. In particular, the case $q = 2$ corresponds to
the scaling predicted for Okuda-Dawson diffusion; see Okuda-Dawson \cite{Okuda1973} and Drake-Greenwood-Navratil-Post \cite{Drake1977}. The Sobolev critical case $q = q_s:=\frac{n+2}{n-2}$ ($1+q=2^*=\frac{2n}{n-2}$) in dimension $n \ge 3$ describes the evolution of a conformal metric by the Yamabe flow. Traditionally, one considers the fast diffusion equation with nonnegative initial condition
\be \label{eq:initial fde}
u(x,0)=u_0(x) \ge 0  \quad   \mbox{on }  \Omega \times \{ t=0\}.
\ee
\eqref{eq:fde}-\eqref{eq:initial fde} is the classical fast diffusion problem on bounded domains and the corresponding nonnegative solutions have already been well studied.

Mass conservation is not allowed under the Dirichlet boundary condition for any $q \in (1,+\infty)$, in contrast with what happens for the Cauchy problem in the whole space. In the Dirichlet case, it is known from the work of Sabinina \cite{Sabinina1,Sabinina2} that solutions of \eqref{eq:fde} with positive initial data extinguish after a finite time $T>0$. For the Cauchy problem in dimensions $n\ge3$, mass is preserved in the range $1<q\le q_c:=n/(n-2)$, whereas finite-time extinction can occur in the very fast diffusion range $q>q_c$. In dimensions $n=1,2$, there is no finite threshold of this form. Extinction-time estimates are discussed by B\'enilan and Crandall \cite{Benilan}.

Positivity of the scalar solution in the Sobolev subcritical range ($1<q<q_s=(n+2)/(n-2)$ for $n\ge3$, and every finite $q>1$ for $n=1,2$) is quantified by a global Harnack inequality:

For all $t_0>0 $ there exists a constant $C>0$, depending on $n,q,\om,t_0$ and $u_0$, such that
\[
\frac{1}{C}d(x)\left(  T-t   \right)^{\frac{1}{q-1}} \le u(x,t) \le Cd(x)\left(  T-t   \right)^{\frac{1}{q-1}}      \quad \mbox{for all } x \in \overline{\om} , t \in[t_0 , T]
\]
where $T>0$ is the extinction time and $d(x)=\dist(x,\pa \om)$, see Kwong \cite{kwong1988,kwong1990} and DiBenedetto-Kwong-Vespri \cite{DKV1991}. This result shows that solutions of the fast diffusion problem have infinite speed of propagation, which means that if the initial data $u_0 \in C_0^1(\om)$ is nonnegative and not identically zero, then the solution will be positive everywhere in $\Omega$ immediately.

Unlike the heat equation ($q=1$), whose nontrivial nonnegative solutions do not extinguish in finite time, the fast diffusion equation ($q>1$) exhibits finite-time extinction. In contrast, the porous medium equation ($0<q<1$ in the present variables) is characterized by finite speed of propagation. For the Dirichlet fast diffusion problem \eqref{eq:fde}--\eqref{eq:initial fde}, we study the asymptotic behavior near extinction through separable solutions. We define the finite
extinction time of a solution $u$ with initial datum $u_0$ as 
\[
T = T(u_0) := \inf  \{t_1 > 0 : u(x,t) = 0 \,  \mbox{ for      a.e. }   x \in \om ,  t \ge t_1\} .
\] 
If $T>0$ is the finite extinction time and $S(x)$ is a solution to the celebrated semilinear elliptic equation of Lane-Emden-Fowler \cite{fowler,lane}
\be \label{eq:steady equation of LEF}
\left\{
\begin{aligned}
 -\Delta S & = cS^q \quad                  \mbox{in }   \om ,\\
S &=0  \quad  \mbox{on } \pa \om  ,\\
S &>0  \quad  \mbox{in }  \om 
\end{aligned}
\right.
\ee
with $c=c(q)=\frac{q}{q-1}$, then
\[
(T-t)^{\frac{1}{q-1}}S(x)
\]
is a separable solution to \eqref{eq:fde}. 

The stability of separable solutions is one of the central problems in fast diffusion on bounded domains. The experiments of Drake, Greenwood, Navratil and Post \cite{Drake1977} on hydrogen plasma diffusion across an octupole magnetic field give a useful picture. After a short transient, the density takes on a fixed spatial shape and then decays towards extinction \cite{2024Bonforte}. The shape stays, while its amplitude falls. This is the behavior described by the scalar fast diffusion model with $q=2$ in the present variables; see also Okuda and Dawson \cite{Okuda1973}.

This observation motivates the following question for the coupled system.
\begin{center}
\small\itshape When both components share a common finite extinction time, what spatial profiles can they approach?
\end{center}
\noindent For two interacting components, the corresponding picture would have two spatial profiles and a common decaying amplitude. Under the common extinction and regularity hypotheses below, we obtain a pair of profiles along a sequence of rescaled times, even when the initial data are not separable. Both profiles are nonnegative, and neither is identically zero. Whether both must be strictly positive remains a further question. Sharing a time factor does not, by itself, make the spatial profiles equal. There is, however, a simple restriction: if the stationary profiles are proportional, they must be equal. Unequal scaled copies of the same shape are therefore ruled out.

For the scalar equation, this stability is studied through the following rescaling. Let
\[
w(x,s)=u(x,t)/(T-t)^{\frac{1}{q-1}} \ \mbox{ and }\  t=T(1-e^{-s}), 
\]
then
\be\label{eq:fast diffusion equation scale}
  \left\{
\begin{aligned}
 \pa_s w^q & =\Delta w + c w^q \quad                  \mbox{in }   \om \times (0,\infty),\\
w &=0  \quad  \mbox{on } \pa \om \times (0,\infty) ,\\
w(0)  &    = u_0/(T)^{\frac{1}{q-1}}  \quad  \mbox{on }  \om \times \{ s=0 \} .
\end{aligned}
\right.
\ee
In 1980, Berryman and Holland \cite{Berryman} proved that, in the Sobolev subcritical range and under a regularity assumption, the rescaled scalar solution satisfies
\[
\|w(s_j)-S\|_{H_0^1(\Omega)}\longrightarrow0
\quad\text{for some sequence }s_j\to\infty.
\]
In 2000, Feireisl and Simondon \cite{Eduard2000Convergence} established uniform convergence. In 2012, Bonforte, Grillo and V\'azquez \cite{Grillo0Behaviour} proved convergence of the relative error. Jin and Xiong resolved the regularity issue raised by Berryman and Holland through their work on optimal Dirichlet boundary regularity, published in 2023 and 2025 \cite{2019Optimal,JX22}; their 2023 paper \cite{JX20} also studied extinction profiles and algebraic convergence rates. Sharp rates were obtained by Bonforte and Figalli in 2021 \cite{2019Sharp} and Akagi in 2023 \cite{2021Sharprate}. In 2023, Choi, McCann and Seis \cite{CMS} developed finer asymptotics.

For a broader account of the scalar problem, we refer to the survey of Bonforte and Figalli \cite{2024Bonforte}. Two properties will help connect the scalar equation with our system. The celebrated B\'enilan--Crandall estimates \cite{Benilan} for solutions to the fast diffusion equation \eqref{eq:fde} with $q\in (1,+\infty)$ hold in the distributional sense and read
\be \label{eq:Benilan-Crandall-fde}
 \pa_t u \le \frac{u}{(q-1)t} .
\ee
This is just the weak formulation of the following time monotonicity:
 \[
 t \to  t^{-\frac{1}{q-1}} u(x,t) \quad
 \mbox{is nonincreasing for a.e.} \, x \in \om.
 \]
This monotonicity plays a key role in the regularity estimates, see \cite{Benilan,Vaz2007} for more details. 

Next, we recall a monotonicity property of the nonlinear Rayleigh quotient, first proved by Berryman-Holland
in their pioneering paper \cite{Berryman}. Consider a nonlinear Rayleigh
quotient of \eqref{eq:fde} for all $q \in (1, +\infty)$
\[
Q[f] := \frac{\| \nabla f  \|^2_{L^2(\om)}}{\| f  \|^2_{L^{q +1}(\om)}}.
\]
It has been shown that $t \mapsto Q[u(t)]$ is nonincreasing along the flow of \eqref{eq:fde}.

In the scalar problem, this monotonicity gives upper and lower bounds near extinction. For a system, looking at one component at a time misses the joint effect of the two loss terms.

The same idea of a fixed profile appears in boundary diffusion. In joint work with Jin and Xiong \cite{2024JXY}, published in 2024, we studied a nonlinear boundary diffusion problem in which $u$ is harmonic in $\Omega$ and evolves on the boundary according to $\partial_t(u^r)=-\partial_\nu u-bu$. When $n\ge3$, at the critical exponent $r=n/(n-2)$ and with $b=-(n-2)H/[2(n-1)]$, where $H$ is the boundary mean curvature, this equation is the unnormalized boundary Yamabe flow; see Brendle \cite{brendle2002generalization}. Our positive-solution bounds and characterization of the asymptotic behavior by separable profiles, including sharp convergence rates, were established in the Sobolev subcritical range. In that problem, the evolution takes place at the boundary. Here the interaction is in the interior, and we ask what it does to the profiles near extinction.

In this paper, we study the stability of separable solutions to system \eqref{eq:fast diffusion system} in analogy with the scalar fast diffusion equation \eqref{eq:fde}, in the range of Sobolev subcritical  exponents, i.e.,
\begin{equation}\label{eq:sobolev subcritical exponent}
\left\{
  \begin{aligned}
  & p\in  (1,+\infty)  \quad  &\mbox{if }   n = 1,2, \\
  &  p\in  \left(1,  \frac{n+2}{n-2}\right) \quad   &\mbox{if }    n\geq3.
\end{aligned}
\right.
 \end{equation}
For a separable pair $(u,v)=(b\Phi,b\phi)$, the common time factor is
\[
b(t)=(T_0-t)^{1/(p-1)},\qquad 0\le t<T_0,
\]
and it solves the ODE
\[
\frac{d}{dt}[b(t)]^p+\frac{p}{p-1}b(t)=0
\]
with the initial data $b(0)=T_0^{\frac{1}{p-1}}$ for some $T_0>0$, and $(\Phi,\phi)$ is a pair of nonnegative solutions of the nonlinear stationary system with Dirichlet boundary condition
\be \label{eq:steady system}
\left\{
\begin{aligned}
 -\Delta \Phi + \alpha \phi & = \frac{p}{p-1}\Phi^p \quad                  \mbox{in }   \om ,\\
 -\Delta \phi + \alpha  \Phi & = \frac{p}{p-1}\phi^p  \quad                  \mbox{in }   \om ,\\
 \Phi =\phi&=0  \quad  \mbox{on } \pa \om  .
\end{aligned}
\right.
\ee
The convergence theorem below yields a nonnegative classical solution pair $(\Phi,\phi)$ of \eqref{eq:steady system}, with neither component identically zero.

The loss terms still allow scalar upper comparison, and the B\'enilan--Crandall estimate holds on the diagonal $u=v$. Away from the diagonal, the comparison argument runs into a sign problem: differences of the coupling terms have no fixed sign. We therefore need a common dissipation structure that keeps track of both components.

For the nonnegative solutions $(u,v)$ to \eqref{eq:fast diffusion system}, we introduce the corresponding Rayleigh quotient:
\[
Q_\alpha [f,g]=  \frac{  \displaystyle \int_{\Omega}|\nabla f|^{2}\,\ud x    +\int_{\Omega}|\nabla g|^{2}\,\ud x   +\int_{\Omega}2\alpha  fg\,\ud x       }{ \displaystyle\left(\int_{ \Omega}f^{p+1} \,\ud x  + \int_{ \Omega}g^{p+1} \,\ud x\right)^{\frac{2}{p+1}}   }  .
\]
The cross term $2\alpha fg$ brings the two losses into one energy. The resulting Rayleigh quotient is nonincreasing and gives upper and lower bounds near the common extinction time. Under the hypotheses below, a Lyapunov functional then gives subsequential strong convergence of the rescaled pair in $H_0^1(\Omega)\times H_0^1(\Omega)$. Strict positivity may be lost in the limit, but the two components cannot both disappear: the integral lower bound prevents this. Nor can just one vanish identically, because the coupling has $\alpha>0$. We also construct a positive diagonal stationary pair by constrained minimization.

An interior zero and the extinction of a whole component are different events. To keep them separate, we use the following times. For continuous solution representatives, define
\[
\begin{aligned}
\sigma_u&:=\inf\{t>0:\ u(x,t)=0\text{ for some }x\in\Omega\},\\
\sigma_v&:=\inf\{t>0:\ v(x,t)=0\text{ for some }x\in\Omega\},
\end{aligned}
\]
with $\inf\varnothing=+\infty$, on the solution branch under consideration, including its terminal trace when available. The joint-positivity time is
\[
T_1:=\min\{\sigma_u,\sigma_v\}.
\]
We retain $T_u,T_v$ for extinction of the entire components and $T_2$ for complete extinction:
\[
T_u:=\inf\{t>0:\ \|u(t)\|_{L^{p+1}(\Omega)}=0\},\qquad
T_v:=\inf\{t>0:\ \|v(t)\|_{L^{p+1}(\Omega)}=0\},
\]
\[
T_2:=\inf\{t>0:\ \|(u(t),v(t))\|_{L^{p+1}(\Omega)}=0\}.
\]
At terminal times, the norms are understood as limits along the solution branch. Whenever these times are defined,
\[
T_1\le T_u\le T_2,\qquad T_1\le T_v\le T_2.
\]

The time estimates in Section~\ref{sec:estimates of time} give, for $p>1$,
\[
T_2 \ge c_{p,\alpha}^{-1}\| (u_0,v_0) \|^{p-1}_{L^p_{\Phi_1}(\om)}
\]
where
\[
\| (u,v) \|_{L^p_{\Phi_1}(\om)} = \left(   \int_{\om} \left( u^p+v^p  \right)\Phi_1 \, \ud x     \right)^{\frac{1}{p}}
\]
and
\[
c_{p,\alpha} :=(p-1)(\lambda_1 +\alpha ) \| \Phi_1  \|^{\frac{p-1}{p}}_{L^1(\om)} ,
\]
see \eqref{eq:lowerbound T_2 Lp_Phi}. $\Phi_1(x)$ denotes the first eigenfunction of the Dirichlet Laplacian on $\om$ and $\lambda_1$ the corresponding eigenvalue, see \eqref{eq:Phi_1 eq}.
For $l>p>1$ and $ l \ge l_c:=\frac{n}{2}(p-1)$
\[
T_1 \le c_{l,p} \left( \min \{\| u_0 \|_{L^l(\om)},\| v_0 \|_{L^l(\om)}   \}  \right)   ^{p-1},
\]
where 
\[
c_{l,p}:=\frac{p(l-p+1)^2}{4(p-1)(l-p)}\mathcal{S}_2^{2}  ,
\]
with $\mathcal{S}_2$ denoting the constant in Sobolev-Poincar\'e inequality \eqref{eq:sobolev-poincare}, see Lemmas \ref{lem:Lp Phi_1 estimates} and \ref{lem:L^l estimates}. Under the common extinction hypothesis stated below, the corresponding lower bounds are
\[
T_u    \ge    \frac{ |  u_0  |^p_{L^p_{\Phi_1}(\om)}  }{  \lambda_1 |  u_0  |_{L^p_{\Phi_1}(\om)}  + \alpha |  v_0  |_{L^p_{\Phi_1}(\om)}         }
\]
and
\[
T_v    \ge    \frac{ |  v_0  |^p_{L^p_{\Phi_1}(\om)}  }{  \lambda_1 |  v_0  |_{L^p_{\Phi_1}(\om)}  + \alpha |  u_0  |_{L^p_{\Phi_1}(\om)}         }
\]
where
\[
|  f  |_{L^p_{\Phi_1}(\om)}:= \frac{\| f  \|_{L^p_{\Phi_1}(\om)}}{\| \Phi_1  \|^{\frac{1}{p}}_{L^1(\om)} }= \left(  \frac{\int_{\om} f^p \Phi_1 \, \ud x}{\int_{\om}  \Phi_1 \, \ud x}  \right)^{\frac{1}{p}},
\]
is the $L^p$ semi-norm of average integral with respect to the weight $  \Phi_1 \, \ud x $,
see Corollary \ref{cor:lower bound of Lp-short time} and Remark \ref{rem:lower bound of T_u and T_v}.

We focus on solution pairs whose positive phase ends in complete extinction, so that $T_1=T_2$. The preceding estimates give the necessary condition
\[
c_{l,p} \left( \min \{\| u_0 \|_{L^l(\om)},\| v_0 \|_{L^l(\om)}   \}  \right)   ^{p-1} \ge c_{p,\alpha
}^{-1}\| (u_0,v_0) \|^{p-1}_{L^p_{\Phi_1}(\om)},
\]
which provides a condition on the initial data $(u_0,v_0)$ as
\be \label{eq:necessary condition on initial data}
\frac{   | (u_0,v_0) |_{L^p_{\Phi_1}(\om)}     }{\min \left\{ \| u_0 \|_{L^l(\om)}   ,\| v_0 \|_{L^l(\om)}    \right\}} \le \left(   \frac{p(\lambda_1 +\alpha)(l-p+1)^2}{4(l-p)}   \mathcal{S}_2^{2} \right)^{\frac{1}{p-1}} 
\ee
for $l>p>1$ and $ l \ge l_c:=\frac{n}{2}(p-1)$, where $\mathcal{S}_2$ is the Sobolev-Poincar\'e constant given in \eqref{eq:sobolev-poincare} and 
\[
| (u_0,v_0) |_{L^p_{\Phi_1}(\om)}:= \frac{\| (u_0,v_0) \|_{L^p_{\Phi_1}(\om)}}{ \| \Phi_1 \|^{\frac{1}{p}}_{L^1(\om)}} = \left(  \frac{\int_{\om} \left( u_0^p+v_0^p   \right)\Phi_1 \, \ud x}{\int_{\om}  \Phi_1 \, \ud x}  \right)^{\frac{1}{p}}.
\]

For the scalar Dirichlet equation, positivity before extinction is already known. In the Sobolev subcritical range, the lower Harnack bound recalled above says more: away from the initial time, the solution is bounded below by a positive multiple of $d(x)(T-t)^{1/(q-1)}$ \cite{kwong1988,kwong1990,DKV1991,2024Bonforte}. The hypothesis below asks for the same positive phase in the coupled problem. Proving similar lower bounds for the system is a natural next question; here we study the extinction profiles within this positive phase.

\textbf{Common extinction hypothesis.}

We assume that both components are strictly positive in $\Omega$ for $0<t<T^*<\infty$, vanish on $\partial\Omega$, and tend to zero in $L^{p+1}(\Omega)$ as $t\uparrow T^*$. With the terminal profiles set equal to zero, we write
\be \label{eq:extinction surmise}
T^*:=T_1=\sigma_u=\sigma_v=T_u=T_v=T_2.
\ee
The two components may start from different data and have different spatial profiles. Their initial data must satisfy \eqref{eq:necessary condition on initial data}. Rescaling both components with $T^*$ gives the following estimates and stability result.

\begin{thm}\label{thm:integralbounds}
Let $n\ge1$, let $\Omega\subset\R^n$ be a bounded smooth domain, and let $p$ satisfy \eqref{eq:sobolev subcritical exponent}. Let $(u,v)$ be a pair of nonnegative weak solutions to \eqref{eq:fast diffusion system} on $\overline\Omega\times[0,T^*)$ satisfying the common extinction hypothesis \eqref{eq:extinction surmise}. Suppose that $\pa_t u, \nabla u, \nabla \pa_t u, \nabla^2 u,\pa_t v, \nabla v, \nabla \pa_t v, \nabla^2 v$ are of class $C(\overline{\om }\times (0,T^*))$. Assume also that $(u(t),v(t))\to(u_0,v_0)$ strongly in $H_0^1(\Omega)\times H_0^1(\Omega)$ as $t\downarrow0$. Then there exists $C>0$ depending only on $n,p,\Omega$ and $Q_\alpha[u_0,v_0]$ such that
\be  \label{eq:integral bounds}
\frac{1}{C} (T^*-t)^{\frac{1}{p-1}} \le \|(u,v)\|_{L^{p+1}(\om)} \le   C (T^*-t)^{\frac{1}{p-1}}, \quad t \in [0,T^*)
\ee
where
\[
 \|(u,v)\|_{L^{p+1}(\om)} =  \left(   \int_{\om} |u|^{p+1}\,\ud x+           \int_{\om} |v|^{p+1}\,\ud x \right)^{\frac{1}{p+1}}.
\]
\end{thm}

To see the spatial profiles as extinction approaches, we rescale the solution pair in Theorem~\ref{thm:integralbounds} by
\begin{equation} \label{eq:changeofvariable}
 \left\{
\begin{array}{lll}
U(x,\tau)&=u(x,t)/(T^*-t)^{\frac{1}{p-1}} \ \mbox{ and }\  t=T^*(1-e^{-\tau}), \\
V(x,\tau)&=v(x,t)/(T^*-t)^{\frac{1}{p-1}} \ \mbox{ and }\  t=T^*(1-e^{-\tau}). \quad
\end{array}
\right.
\end{equation}
Then $(U,V)$ is a pair of solutions to
\be\label{eq:fast diffusion system scale}
  \left\{
\begin{aligned}
 \pa_\tau U^p & =\Delta U -\alpha V+\frac{p}{p-1}U^p \quad                  \mbox{in }   \om \times (0,\infty),\\
 \pa_\tau V^p  &    =\Delta V -\alpha U+\frac{p}{p-1}V^p\quad                  \mbox{in }   \om \times (0,\infty),\\
 U &=V=0  \quad  \mbox{on } \pa \om \times (0,\infty) ,\\
 \left( U(0),V(0) \right) &    = (u_0/(T^*)^{\frac{1}{p-1}},v_0/(T^*)^{\frac{1}{p-1}})   \quad  \mbox{on }  \om \times \{ \tau=0 \} .
\end{aligned}
\right.
\ee
So, it follows from Theorem \ref{thm:integralbounds} that
\be  \label{eq:integral bounds rescaled}
1/C\le \|(U,V)\|_{L^{p+1}(\om)} \le C, \quad \tau \in [0,\infty),
\ee
for  some $C\ge 1$ which depends only on $n,p,\Omega$ and $Q_\alpha[u_0,v_0]$.

\begin{thm}\label{thm:convergence}
Let $(u,v)$ satisfy the hypotheses of Theorem~\ref{thm:integralbounds}, and let $(U,V)$ be defined by \eqref{eq:changeofvariable}. Then there exists an increasing sequence of times $\{\tau_n\}_{n=1}^{\infty}$ with
\[
\tau_n \to \infty    \quad \mbox{as} \quad n \to \infty,
\]
such that
\begin{equation*}
\begin{split}
U(\cdot, \tau_n) \to \Phi(\cdot)  \quad \mbox{in} \quad H_0^{1}{(\om)}, \\
  V(\cdot, \tau_n) \to \phi(\cdot) \quad \mbox{in} \quad H_0^{1}{(\om)}
\end{split}
\end{equation*}
for $\tau_n \to \infty$, where $(\Phi,\phi)$ is a nonnegative classical solution pair of \eqref{eq:steady system}, with $\Phi\not\equiv0$ and $\phi\not\equiv0$.

\end{thm}

Section \ref{sec:comparison principle} establishes scalar upper comparison and the B\'enilan--Crandall estimate on the diagonal. Sections \ref{sec:estimates of time} and \ref{sec:ex} discuss time estimates, local bounds, and evolution and stationary solutions on the diagonal. Sections \ref{sec:integral} and \ref{sec:con} establish the upper and lower bounds near extinction and the subsequential convergence.

\section{\texorpdfstring{Scalar comparison and B\'enilan--Crandall estimates}{Scalar comparison and Benilan--Crandall estimates}}\label{sec:comparison principle}

We work on an interval $(0,T)$ where the components are nonnegative. Each component is bounded above by a scalar fast diffusion solution. On the diagonal $u=v$, the time scaling also gives the scalar B\'enilan--Crandall estimate.

\bpr \label{prop:comparison}
Let $(u,v)$ be a nonnegative solution pair of \eqref{eq:fast diffusion system} on $(0,T)$, with sufficient regularity for scalar Dirichlet comparison. Let $\bar u,\bar v$ solve the two uncoupled scalar fast diffusion problems
\[
\left\{\begin{aligned}
\partial_t(\bar u^p)&=\Delta\bar u,\qquad
\partial_t(\bar v^p)=\Delta\bar v &&\text{in }\Omega\times(0,\infty),\\
\bar u&=\bar v=0 &&\text{on }\partial\Omega\times(0,\infty),\\
(\bar u,\bar v)(x,0)&=(u_0,v_0)(x) &&\text{for }x\in\Omega.
\end{aligned}\right.
\]
Then $u\le\bar u$ and $v\le\bar v$ on $\Omega\times(0,T)$. If both components are strictly positive there, then $T\le\min\{\bar T_u,\bar T_v\}$, where $\bar T_u,\bar T_v$ are the scalar extinction times.
\epr
\bpf
Subtract the equation for $\bar u$ from the first equation of \eqref{eq:fast diffusion system}:
\[
\partial_t(u^p-\bar u^p)-\Delta(u-\bar u)=-\alpha v\le0.
\]
The two differences have the same sign since $s\mapsto s^p$ is increasing on $[0,\infty)$. Writing $z_+=\max\{z,0\}$, the scalar positive-part inequality (obtained by testing with a smooth approximation of $\mathbf{1}_{\{u>\bar u\}}$) yields
\[
\partial_t (u^p-\bar u^p)_+\le\Delta(u-\bar u)_+.
\]
Test with the positive first Dirichlet eigenfunction $\Phi_1$, with eigenvalue $\lambda_1$. The zero boundary data give
\[
\frac{\ud}{\ud t}\int_\Omega (u^p-\bar u^p)_+\Phi_1\,\ud x
\le -\lambda_1\int_\Omega (u-\bar u)_+\Phi_1\,\ud x\le0.
\]
The integral is initially zero, so $u\le\bar u$. The second equation gives $v\le\bar v$ in the same way. If, for instance, $\bar T_u<T$, then $0\le u\le\bar u=0$ for $\bar T_u<t<T$, contradicting strict positivity. The same argument gives $T\le\bar T_v$.
\epf

\bpr \label{prop:benilan-crandall estimates}
Let $(w,w)$ be a smooth positive solution pair of \eqref{eq:fast diffusion system} on $(0,T)$. Then
\begin{equation}\label{eq:benilan-crandall u and v}
\partial_t w(x,t)\le \frac{w(x,t)}{(p-1)t}\qquad (x\in\Omega,\ 0<t<T).
\end{equation}
In particular, $t^{-1/(p-1)}w(x,t)$ is nonincreasing in $t$ for each $x\in\Omega$.
\epr
\bpf
On the diagonal, $\partial_t w^p=\Delta w-\alpha w$. Set $a=1/(p-1)$ and $w_\lambda(x,t)=\lambda^{-a}w(x,\lambda t)$ for $\lambda>1$ and $0<t<T/\lambda$. Since $1-ap=-a$,
\[
\begin{aligned}
\partial_t w_\lambda^p(x,t)
&=\lambda^{1-ap}(\partial_t w^p)(x,\lambda t)\\
&=\lambda^{-a}(\Delta w-\alpha w)(x,\lambda t)
=\Delta w_\lambda(x,t)-\alpha w_\lambda(x,t).
\end{aligned}
\]
Also $w_\lambda(x,0)=\lambda^{-a}w_0(x)\le w_0(x)$, and both functions have zero boundary data. Subtracting the scalar equations gives
\[
\partial_t(w_\lambda^p-w^p)-\Delta(w_\lambda-w)=-\alpha(w_\lambda-w).
\]
On $\{w_\lambda>w\}$ the right-hand side is nonpositive, so the positive-part argument above gives $w_\lambda\le w$. For fixed $0<t<T$, take $h>0$ with $t+h<T$ and set $\lambda=1+h/t$. Then
\[
\frac{w(x,t+h)-w(x,t)}{h}
\le \frac{(1+h/t)^a-1}{h}\,w(x,t)
\ \longrightarrow\ \frac{a}{t}w(x,t)\qquad (h\downarrow0).
\]
This proves \eqref{eq:benilan-crandall u and v}. The asserted monotonicity follows by differentiating $t^{-a}w(x,t)$.
\epf

\brem
For general smooth positive solution pairs, does
\[
\min\left\{\partial_t u-\frac{u}{(p-1)t},\
                    \partial_t v-\frac{v}{(p-1)t}\right\}\le0
\]
hold at every point of the common positivity interval? Proposition~\ref{prop:benilan-crandall estimates} answers this question on the diagonal, where both quantities are nonpositive. For general pairs, the pointwise bound is not established here. The component satisfying the bound may change with the space--time point, so even this bound would not by itself imply that either $t^{-1/(p-1)}u$ or $t^{-1/(p-1)}v$ is nonincreasing throughout a time interval.
\erem

\section{\texorpdfstring{Estimates of the joint-positivity time $T_1$ and complete extinction time $T_2$}{Estimates of the joint-positivity time T1 and complete extinction time T2}} \label{sec:estimates of time}
All estimates below are understood on intervals where the solution pair exists and is nonnegative. Statements at an extinction time require the corresponding terminal norm limit.
Let $\Phi_1(x)$ denote the first eigenfunction of the Dirichlet Laplacian on $\om$, satisfying
\begin{equation}\label{eq:Phi_1 eq}
 \left\{
\begin{aligned}
 -\Delta  \Phi_1 &=\lambda_1 \Phi_1 \quad                           \mbox{in }   \Omega ,\\
           \Phi_1 & =  0  \quad   \mbox{on } \partial \Omega
\end{aligned}
\right.
\end{equation}
where $\lambda_1>0$ denotes the corresponding first eigenvalue. The eigenfunction $\Phi_1\in C^{\infty}(\overline{\om})$ is positive in $\om$ and satisfies
\be \label{eq:positiveity property of Phi_1}
\Phi_1(x) \asymp \dist(x,\pa \om),
\ee
where
\[
 f \asymp g
\]
means there is a constant $C>1$ such that
\[
\frac{1}{C}f \le g \le C f.
\]
Hence, it is convenient to use it as a smooth extension of the distance function. We normalize it by
\[
\| \Phi_1(x) \|_{L^2(\om)} =1
\]
without loss of generality.

\blem\label{lem:Lp Phi_1 estimates}($L^p_{\Phi_1} $-Estimates)

Let $(u,v)$ be a pair of nonnegative weak solutions to \eqref{eq:fast diffusion system} with $(u_0,v_0) \in L^p_{\Phi_1}(\om)$ and $p \in (1,+\infty)$. Then the complete extinction time $T_2=T_2(u_0,v_0)$ can be estimated from below as
\be \label{eq:lowerbound T_2 Lp_Phi}
T_2 \ge c_{p,\alpha}^{-1}\| (u_0,v_0) \|^{p-1}_{L^p_{\Phi_1}(\om)}  \quad   \mbox{ where }  \quad c_{p,\alpha} :=(p-1)(\lambda_1 +\alpha ) \| \Phi_1  \|^{\frac{p-1}{p}}_{L^1(\om)} ,
\ee
\[
\| (u,v) \|_{L^p_{\Phi_1}(\om)} = \left(   \int_{\om} \left( u^p+v^p  \right)\Phi_1 \, \ud x     \right)^{\frac{1}{p}}
\]
Moreover, for all $0 \le t_0 \le t \le T_2$,
\be \label{eq:Lp Phi_1 estimates}
\left(  \int_{\om} \left(u(t_0)^p+ v(t_0)^p   \right)\Phi_1  \, \ud x \right)^{\frac{p-1}{p}} - c_{p,\alpha}(t-t_0) \le \left(  \int_{\om} \left(u(t)^p+ v(t)^p   \right)\Phi_1  \, \ud x \right)^{\frac{p-1}{p}} \le c_{p,\alpha}(T_2 -t).
\ee

\elem

\bpf
We differentiate in time to get
\[
\frac{\ud}{\ud t} \int_{\om} u^p \Phi_1 \, \ud x  = \int_{\om} \left( u \Delta \Phi_1 - \alpha v\Phi_1  \right) \, \ud x=\int_{\om} \left( -\lambda_1 u  \Phi_1 - \alpha v\Phi_1  \right) \, \ud x,
\]
\[
\frac{\ud}{\ud t} \int_{\om} v^p \Phi_1 \, \ud x  = \int_{\om} \left( v \Delta \Phi_1 - \alpha u\Phi_1  \right) \, \ud x=\int_{\om} \left( -\lambda_1 v  \Phi_1 - \alpha u\Phi_1  \right) \, \ud x.
\]
So
\[
\frac{\ud}{\ud t} \int_{\om} \left( u^p +v^p  \right) \Phi_1 \, \ud x  = -(\lambda_1 +\alpha ) \int_{\om} \left( u+v\right) \Phi_1 \, \ud x.
\]
Using H\"older's inequality, one obtains
\begin{equation*}
\begin{split}
\frac{\ud}{\ud t} \int_{\om} \left( u^p +v^p  \right) \Phi_1 \, \ud x  & \ge  -(\lambda_1 +\alpha ) \| \Phi_1  \|^{\frac{p-1}{p}}_{L^1{(\om)}}\left(    \left(\int_{\om} u^p \Phi_1 \, \ud x     \right)^{\frac{1}{p}} +\left(\int_{\om} v^p \Phi_1 \, \ud x     \right)^{\frac{1}{p}}\right)\\
& \ge -p(\lambda_1 +\alpha ) \| \Phi_1  \|^{\frac{p-1}{p}}_{L^1{(\om)}}  \left(    \int_{\om} u^p \Phi_1 \, \ud x      +\int_{\om} v^p \Phi_1 \, \ud x   \right)^{\frac{1}{p}}.
\end{split}
\end{equation*}

Integrating over $[t_0,t]$ gives the first inequality in \eqref{eq:Lp Phi_1 estimates}. Applying this bound on $[t,T_2]$ and using the vanishing of both components at $T_2$ gives the second inequality. Taking $t_0=0$ and $t=T_2$ gives \eqref{eq:lowerbound T_2 Lp_Phi}.

\epf

\begin{cor}\label{cor:lower bound of Lp-short time}
Let $(u,v)$ be a solution pair as in Lemma~\ref{lem:Lp Phi_1 estimates}, with nonnegative initial data $u_0,v_0$, neither of which is identically zero. Then there exists a time $t_0=t_0(\alpha,\om,p,u_0,v_0)>0$ such that
\be \label{eq:short time lower bound of Lp}
\min\left\{ \int_{\om} u^p  \, \ud x , \int_{\om} v^p  \, \ud x \right\} \ge C >0 , \quad \mbox{for all } t \in [0,t_0].
\ee

\end{cor}

\bpf
From the argument in Lemma \ref{lem:Lp Phi_1 estimates}, we know that
\[
\int_{\om} u^p\Phi_1 \, \ud x    \quad \mbox{and} \quad \int_{\om} v^p\Phi_1 \, \ud x
\]
are both decreasing in time and
\begin{equation*}
\begin{split}
\frac{\ud}{\ud t} \int_{\om} u^p \Phi_1 \, \ud x  & =\int_{\om} \left( -\lambda_1 u  \Phi_1 - \alpha v\Phi_1  \right) \, \ud x  \\
& \ge - \| \Phi_1  \|^{\frac{p-1}{p}}_{L^1{(\om)}}\left(  \lambda_1  \left(\int_{\om} u^p \Phi_1 \, \ud x     \right)^{\frac{1}{p}} +\alpha \left(\int_{\om} v^p \Phi_1 \, \ud x     \right)^{\frac{1}{p}}\right) \\
& \ge - \| \Phi_1  \|^{\frac{p-1}{p}}_{L^1{(\om)}}\left(  \lambda_1  \left(\int_{\om} u_0^p \Phi_1 \, \ud x     \right)^{\frac{1}{p}} +  \alpha  \left(\int_{\om} v_0^p \Phi_1 \, \ud x     \right)^{\frac{1}{p}}\right).
\end{split}
\end{equation*}
Similarly,
\[
\frac{\ud}{\ud t} \int_{\om} v^p \Phi_1 \, \ud x  \ge - \| \Phi_1  \|^{\frac{p-1}{p}}_{L^1{(\om)}}\left(  \lambda_1  \left(\int_{\om} v_0^p \Phi_1 \, \ud x     \right)^{\frac{1}{p}} +\alpha \left(\int_{\om} u_0^p \Phi_1 \, \ud x     \right)^{\frac{1}{p}}\right).
\]
Integrating the above two inequalities from $0$ to $t$, we obtain
\[
\int_{\om} u^p \Phi_1 \, \ud x  \ge \int_{\om} u_0^p \Phi_1 \, \ud x - t \| \Phi_1  \|^{\frac{p-1}{p}}_{L^1{(\om)}}\left(  \lambda_1  \left(\int_{\om} u_0^p \Phi_1 \, \ud x     \right)^{\frac{1}{p}} +\alpha \left(\int_{\om} v_0^p \Phi_1 \, \ud x     \right)^{\frac{1}{p}}\right)
\]
and
\[
\int_{\om} v^p \Phi_1 \, \ud x  \ge \int_{\om} v_0^p \Phi_1 \, \ud x - t \| \Phi_1  \|^{\frac{p-1}{p}}_{L^1{(\om)}}\left(  \lambda_1  \left(\int_{\om} v_0^p \Phi_1 \, \ud x     \right)^{\frac{1}{p}} +\alpha \left(\int_{\om} u_0^p \Phi_1 \, \ud x     \right)^{\frac{1}{p}}\right).
\]
\eqref{eq:short time lower bound of Lp} is proved if we choose proper $t_0(\alpha,\om,p,u_0,v_0)$ small enough and use the property \eqref{eq:positiveity property of Phi_1} of $\Phi_1$.
\epf

\blem\label{lem:L^l estimates}($L^l $-Estimates)

Let $(u,v)$ be a pair of nonnegative weak solutions to \eqref{eq:fast diffusion system} with $(u_0,v_0) \in L^l(\om)$ and $p \in (1,+\infty)$. If $l>p$ and $ l \ge l_c:=\frac{n}{2}(p-1)$, then the extinction times $T_u=T_u(u_0,v_0)$ and $T_v=T_v(u_0,v_0)$ can be estimated from above as
\be \label{eq:upperbound T_u,T_v L^l}
T_u \le c_{l,p}\| u_0 \|^{p-1}_{L^l(\om)}, \ T_v \le c_{l,p}\| v_0 \|^{p-1}_{L^l(\om)}
  \quad   \mbox{ where }  \quad c_{l,p}:=\frac{p(l-p+1)^2}{4(p-1)(l-p)}\mathcal{S}_2^{2}  ,
\ee
and $\mathcal{S}_2^{2}$ is the Sobolev-Poincar\'e constant given in \eqref{eq:sobolev-poincare}. Moreover, for $0\le t_0\le t\le T_u$,
\be \label{eq:u L^l   estimates}
c_{l,p}^{-1}(T_u -t)  \le \| u(t)  \|^{p-1}_{L^l(\om)} \le  \| u(t_0)  \|^{p-1}_{L^l(\om)}  - c_{l,p}^{-1}(t -t_0) ,
\ee
and, for $0\le t_0\le t\le T_v$,
\be \label{eq:v L^l   estimates}
c_{l,p}^{-1}(T_v -t)  \le \| v(t)  \|^{p-1}_{L^l(\om)} \le  \| v(t_0)  \|^{p-1}_{L^l(\om)}  - c_{l,p}^{-1}(t -t_0) .
\ee

\elem

\bpf
Recall that
\begin{equation}\label{eq:u^l}
\begin{split}
\frac{\ud}{\ud t} \int_{\om} u^l  \, \ud x & = -\frac{4l(l-p)}{p(l-p+1)^2}\int_{\om} |\nabla u ^{\frac{l-p+1}{2}}|^2 \, \ud x-  \frac{\alpha l}{p} \int_{\om} u^{l-p}v \, \ud x  \\
& \le -\frac{4l(l-p)}{p(l-p+1)^2}\int_{\om} |\nabla u ^{\frac{l-p+1}{2}}|^2 \, \ud x \\
& \le -\frac{4l(l-p)}{p(l-p+1)^2} \mathcal{S}_2^{-2} \left( \int_{\om} u^l \, \ud x \right)^{1-\frac{p-1}{l}},
\end{split}
\end{equation}
\begin{equation}\label{eq:v^l}
\begin{split}
\frac{\ud}{\ud t} \int_{\om} v^l  \, \ud x & = -\frac{4l(l-p)}{p(l-p+1)^2}\int_{\om} |\nabla v ^{\frac{l-p+1}{2}}|^2 \, \ud x- \frac{\alpha l}{p} \int_{\om} v^{l-p}u \, \ud x  \\
& \le -\frac{4l(l-p)}{p(l-p+1)^2}\int_{\om} |\nabla v ^{\frac{l-p+1}{2}}|^2 \, \ud x \\
& \le -\frac{4l(l-p)}{p(l-p+1)^2} \mathcal{S}_2^{-2} \left( \int_{\om} v^l \, \ud x \right)^{1-\frac{p-1}{l}}
\end{split}
\end{equation}
where $\mathcal{S}_2$ is the constant in Sobolev-Poincar\'e inequality
\be \label{eq:sobolev-poincare}
\| f  \|_{L^{\frac{2l}{l-p+1}}(\om)} \le \mathcal{S}_2 \| \nabla f \|_{L^2(\om)}  ,  \quad \forall f \in H_0^1(\om).
\ee
For $n\ge3$, notice that
\[
\frac{2l}{l-p+1} \le 2^*:=\frac{2n}{n-2}    \quad \mbox{if and only if}   \quad  l \ge l_c:=\frac{n}{2}(p-1).
\]
Integrating the above inequalities \eqref{eq:u^l} and \eqref{eq:v^l} over $[t_0,t]$, we get
\[
\left( \int_{\om} u(x,t)^l \, \ud x \right)^{\frac{p-1}{l}} \le \left( \int_{\om} u(x,t_0)^l \, \ud x \right)^{\frac{p-1}{l}} - \frac{4(p-1)(l-p)}{p(l-p+1)^2} \mathcal{S}_2^{-2}  (t - t_0),
\]
\[
\left( \int_{\om} v(x,t)^l \, \ud x \right)^{\frac{p-1}{l}} \le \left( \int_{\om} v(x,t_0)^l \, \ud x \right)^{\frac{p-1}{l}} - \frac{4(p-1)(l-p)}{p(l-p+1)^2} \mathcal{S}_2^{-2}  (t - t_0),
\]
which prove the upper bounds in \eqref{eq:u L^l   estimates} and \eqref{eq:v L^l   estimates}. Applying the same bounds on $[t,T_u]$ and $[t,T_v]$, respectively, and using the vanishing of the corresponding component at its extinction time gives the lower bounds. Finally, taking the starting time to be zero gives \eqref{eq:upperbound T_u,T_v L^l}.

\epf

\brem \label{rem:critical fast diffusion exponent}
For $n\ge3$, define \[
p_c:=1+\frac{2}{n-2}=\frac{n}{n-2}, \quad l_c=\frac{n}{2}(p-1),
\]
one notes that
\[
p\in (1,p_c)   \quad \mbox{if and only if}   \quad  l_c \in (0,p).
\]
Therefore, $(1,p_c)$ is called the good fast diffusion parameter range for the exponent $p$ and we say that $p_c=\frac{n}{n-2}$ is the fast diffusion critical exponent of system \eqref{eq:fast diffusion system} like \eqref{eq:fde}. In this range, $l>p$ implies $l>l_c$, so Lemma \ref{lem:L^l estimates} applies. For $n=1,2$, the embeddings $H_0^1(\Omega)\hookrightarrow L^r(\Omega)$ hold for every finite $r$, and $l_c<p$ for every $p>1$. Thus there is no finite upper threshold $p_c$ in these dimensions.
\erem

\brem \label{rem:critical fast diffusion exponent with sobolev}
For $n\ge3$, if we consider energy solutions $u(t),v(t) \in L^{p+1}(\om)$, then it is natural to require that $p+1 > l_c$, to employ the energy estimate of  Lemma \ref{lem:L^l estimates}.
\[
p+1 > l_c  \quad \mbox{if and only if}     \quad  p<p_s:=\frac{n+2}{n-2}.
\]
As above, $l_c$ provides a connection between fast diffusion critical exponent $p_c$ and the Sobolev one $p_s$. And $p_s:=\frac{n+2}{n-2}$ is relevant to the asymptotic behavior near extinction. When $n=1,2$, $p+1>l_c$ holds for every $p>1$, and every finite $p>1$ is Sobolev subcritical, as in \eqref{eq:sobolev subcritical exponent}. 
\erem

\brem \label{rem:lower bound of T_u and T_v}
Assume the common extinction hypothesis \eqref{eq:extinction surmise}. The integral estimates in the proof of Corollary \ref{cor:lower bound of Lp-short time} give the following lower bounds for $T_u=T_v=T^*$:
\be \label{eq:lower bound T_u  }
T_u    \ge  \frac{  \| u_0  \|_{L^p_{\Phi_1}(\om)}^p }{  \| \Phi_1  \|^{\frac{p-1}{p}}_{L^1{(\om)}} \left(    \lambda_1 \| u_0  \|_{L^p_{\Phi_1}(\om)} +\alpha \| v_0  \|_{L^p_{\Phi_1}(\om)}        \right) } =  \frac{ |  u_0  |^p_{L^p_{\Phi_1}(\om)}  }{  \lambda_1 |  u_0  |_{L^p_{\Phi_1}(\om)}  + \alpha  |  v_0  |_{L^p_{\Phi_1}(\om)}         }
\ee
and
\be  \label{eq:lower bound T_v  }
T_v    \ge  \frac{  \| v_0  \|_{L^p_{\Phi_1}(\om)}^p }{  \| \Phi_1  \|^{\frac{p-1}{p}}_{L^1{(\om)}} \left(    \lambda_1 \| v_0  \|_{L^p_{\Phi_1}(\om)} +\alpha \| u_0  \|_{L^p_{\Phi_1}(\om)}        \right) } =  \frac{ |  v_0  |^p_{L^p_{\Phi_1}(\om)}  }{  \lambda_1 |  v_0  |_{L^p_{\Phi_1}(\om)}  + \alpha |  u_0  |_{L^p_{\Phi_1}(\om)}         }
\ee
where
\[
\| u_0  \|_{L^p_{\Phi_1}(\om)}:= \left( \int_{\om} u_0^p \Phi_1 \, \ud x \right)^{\frac{1}{p}}
\]
and
\[
|  u_0  |_{L^p_{\Phi_1}(\om)}:= \frac{\| u_0  \|_{L^p_{\Phi_1}(\om)}}{\| \Phi_1  \|^{\frac{1}{p}}_{L^1(\om)} }= \left(  \frac{\int_{\om} u_0^p \Phi_1 \, \ud x}{\int_{\om}  \Phi_1 \, \ud x}  \right)^{\frac{1}{p}}.
\]
Letting $t\uparrow T^*$ in the weighted integral inequalities gives these bounds. Under the common extinction hypothesis, $\sigma_u=\sigma_v=T^*$, so the bounds also apply to the first interior contact times. Without this hypothesis, the integral estimates alone give no such lower bounds for $\sigma_u$ or $\sigma_v$.

\erem

\section{\texorpdfstring{A priori estimates and existence on the diagonal}{A priori estimates and existence on the diagonal}}\label{sec:ex}

We give local a priori bounds for solution pairs, construct weak solutions for diagonal initial data, and establish existence of a positive diagonal stationary pair.
\subsection{Fast diffusion system}

We begin with the weak formulation on the time interval where the solution pair is defined.

\begin{defn}\label{defn:weak solution}
We say that $(u,v)$ is a pair of weak solutions of \eqref{eq:fast diffusion system} if
\be
\begin{split}
 &u,v \in L^{\infty} \left( Q    \right),  \mbox{ where }  Q=\om \times (t_1,t_2)\\
  &u,v \ge 0,  \\
   &u^p,v^p \in C \left( [t_1,t_2] ; L^1(\om)   \right)
\end{split}
\ee
satisfy
\be
\begin{split}
& \int_{t_1}^{t_2}\int_{\om}  \Big(  u^p \pa_t \varphi +  u  \Delta \varphi - \alpha v \varphi                 \Big) \, \ud x \ud t  = \int_{\om} u(x,t_2)^p\varphi(x,t_2) \, \ud x - \int_{\om} u(x,t_1)^p\varphi(x,t_1)  \, \ud x,\\
 & \int_{t_1}^{t_2}\int_{\om}  \Big(  v^p \pa_t \varrho + v  \Delta \varrho -  \alpha u \varrho                 \Big)  \, \ud x \ud t  = \int_{\om} v(x,t_2)^p\varrho(x,t_2) \, \ud x - \int_{\om} v(x,t_1)^p\varrho(x,t_1) \, \ud x
\end{split}
\ee
for any closed interval $[t_1,t_2]$ contained in the interval of existence and $ \varphi, \varrho \in C^{2,1}(\overline{\om}\times[t_1,t_2])$ which both vanish on $\pa \om \times [t_1,t_2]$.
\end{defn}

\subsubsection{\texorpdfstring{A priori local $L^\infty$ bounds}{A priori local Linfinity bounds}}

We recall the Sobolev inequalities needed for the $L^\infty$ estimate.

For $0<q<+\infty$ if $n\le 2$ and $0<q \le \frac{n+2}{n-2}$ if $n\ge 3$, there exists $C>0$ depending only on $n,q,\om$ such that
\be \label{eq:Sobolev inequality om}
\left( \int_{\om} |f|^{q+1} \, \ud x     \right)^{\frac{2}{q+1}}    \le C \left( \int_{\om} |\nabla f|^{2} \, \ud x   +  \int_{\om} | f|^{2} \, \ud x    \right)             \quad \forall f \in H^1(\om).
\ee
In particular, we have the Sobolev  inequality on the ball with radius $r>0$. For $0<q<+\infty$ if $n\le 2$ and $0<q\le \frac{n+2}{n-2}$ if $n\ge 3$, there exists $C>0$ depending only on $n,q,r$ such that
\be \label{eq:Sobolev inequality Br}
\left( \int_{B_r} |f|^{q+1} \, \ud x     \right)^{\frac{2}{q+1}}    \le C \left( \int_{B_r} |\nabla f|^{2} \, \ud x   +  \int_{B_r} | f|^{2} \, \ud x    \right)             \quad \forall f \in H^1(B_r).
\ee

The proof given here is based on a variant of Moser's \cite{Moser} iteration technique.

We use the following notation
\[
Q^*=\left\{  (x,t)\in \R^{n+1}   : |x|<2, -4<t<0     \right\},
\]
\[
Q=\left\{  (x,t)\in Q^*   : |x|<1, -1<t<0     \right\},
\]
\[
Q_r=\left\{  (x,t)\in  Q^*   : |x|<r, -r^2 <t<0     \right\},
\]
\[
B_r=\{x \in \R^n :|x|<r\}.
\]

\bthm \label{thm:L^infty local normalize}
Assume $n\ge3$ and that $p$ satisfies \eqref{eq:sobolev subcritical exponent}. Let $(u,v)$ be a pair of smooth nonnegative solutions of the equation
\be\label{eq:fast diffusion system normalize without initial data}
  \left\{
\begin{aligned}
 \pa_t u^p & =\Delta u -\alpha v \quad                  \mbox{in }   Q^*,\\
 \pa_t v^p  &    =\Delta v -\alpha u\quad                  \mbox{in }   Q^*\\
\end{aligned}
\right.
\ee
with
\[
\min \{ \|u\|_{L^{p}(B_{1/2})} ,\|v\|_{L^{p}(B_{1/2}) }     \}   \ge \delta
\]
for some $\delta>0$ and all $t\in(-4,0)$. Then
\be  \label{eq: u,v L infty bound}
\|  u \|_{L^{\infty}(Q)} + \|  v \|_{L^{\infty}(Q)} \le C \|(u,v)\|^{\theta}_{L^{p+1}(Q^*)}= C\left(    \int_{Q^*} \left( u^{p+1 } + v^{p+1} \right) \, \ud x \ud t \right)^{\frac{\theta}{p+1}}
\ee
Here $C$ depends on $n,p,\delta$, and $\theta=(p+1)/(p+1-n(p-1)/2)$.
\ethm

\bpf
We give the proof for $n\ge3$ and fixed $\alpha>0$. We use $\beta\ge0$ for the iteration exponent.

First we show the estimate of $u$. For $1/2<r<R\le2$, let $\xi(x),\eta(t)$ be smooth cut-off functions such that
\be\nonumber
  \left\{
\begin{aligned}
 \xi(x) & =0 \quad    \quad               |x|>R,\\
 0 \le \xi(x)& \le 1 \quad        r\le |x| \le R     ,    \\
 \xi(x) & =1 \quad    \quad               |x|<r,
\end{aligned}
\right.
\ee
\be\nonumber
  \left\{
\begin{aligned}
 \eta(t) & =0 \quad     \quad              t < -R^2,\\
 0 \le \eta(t)& \le 1 \quad        -R^2 \le t \le -r^2     ,    \\
 \eta(t) & =1 \quad   \quad                t>-r^2,
\end{aligned}
\right.
\ee
with
\[
\left|    D\xi    \right| \le \frac{C}{R-r} , \quad \left|    D_t\eta    \right| \le \frac{C}{|R-r|^2}.
\]
Multiplying the first equation of \eqref{eq:fast diffusion system normalize without initial data} by $\xi^2\eta^2u^{1+\beta}$ for $\beta\ge0$ and integrating over $Q^{\tau}_R= B_R \times (-R^2,\tau)$ for $\tau \in (-r^2,0)$, we obtain
\[
\iint_{Q^{\tau}_R} \xi^2 \eta^2 u^{1+\beta} \pa_t u^p  = \iint_{Q^{\tau}_R} \xi^2 \eta^2 u^{1+\beta} \Delta u - \alpha \iint_{Q^{\tau}_R} \xi^2 \eta^2 u^{1+\beta} v.
\]
Then, integration by parts gives
\begin{equation*}
\begin{split}
\iint_{Q^{\tau}_R} \xi^2 \eta^2 u^{1+\beta} \pa_t u^p + \iint_{Q^{\tau}_R} \xi^2 \eta^2 \nabla\left( u^{1+\beta} \right) \cdot \nabla u
& \le -2 \iint_{Q^{\tau}_R} \xi \eta^2 u^{1+\beta} \nabla \xi \cdot \nabla u- \alpha \iint_{Q^{\tau}_R} \xi^2 \eta^2 u^{1+\beta} v\\
& \le -2 \iint_{Q^{\tau}_R} \xi \eta^2 u^{1+\beta} \nabla \xi \cdot \nabla u
\end{split}
\end{equation*}
where we used the nonnegativity of $u,v$. On the left hand side, for the first term, we have
\begin{equation*}
\begin{split}
\iint_{Q^{\tau}_R} \xi^2 \eta^2 u^{1+\beta} \pa_t u^p & = \frac{p}{p+1+\beta} \iint_{Q^{\tau}_R} \xi^2 \eta^2 \pa_t  u^{p+1+\beta}  \\
& = \frac{p}{p+1+\beta} \iint_{Q^{\tau}_R} \pa_t \left( \xi^2 \eta^2   u^{p+1+\beta} \right) - \frac{2p}{p+1+\beta}\iint_{Q^{\tau}_R} \eta D_t \eta \xi^2 u^{p+1+\beta} \\
& \ge \frac{p}{p+1+\beta} \int_{B_R}  \xi^2   u^{p+1+\beta}(\tau)  - \frac{C}{|R-r|^2}\iint_{Q^\tau_R}  u^{p+1+\beta}.
\end{split}
\end{equation*}
For the second term, we have
\begin{equation*}
\begin{split}
\iint_{Q^{\tau}_R} \xi^2 \eta^2 \nabla\left( u^{1+\beta} \right) \cdot \nabla u & = \left(  1+\beta \right) \iint_{Q^{\tau}_R} \xi^2 \eta^2 u^\beta |\nabla u|^2 \\
& = \frac{4(1+\beta)}{(2+\beta)^2}\iint_{Q^{\tau}_R} \xi^2 \eta^2  |\nabla (u^{\frac{2+\beta}{2}})|^2 .
\end{split}
\end{equation*}
On the right hand side, we have
\begin{equation*}
\begin{split}
-2 \iint_{Q^{\tau}_R} \xi \eta^2 u^{1+\beta} \nabla \xi \cdot \nabla u & = -  \iint_{Q^{\tau}_R} \left(   \xi \eta u^{\frac{\beta}{2}} \nabla u    \right)  \cdot \left(   2\eta u^{\frac{2+\beta}{2}}   \nabla \xi  \right) \\
& \le \frac{1+\beta}{2}\iint_{Q^{\tau}_R} \xi^2 \eta^2 u^\beta |\nabla u|^2    +  \frac{1}{2(1+\beta)}  \iint_{Q^{\tau}_R} 4 \eta^2 u^{2+\beta} |\nabla \xi|^2 \\
& = \frac{2(1+\beta)}{(2+\beta)^2} \iint_{Q^{\tau}_R} \xi^2 \eta^2 |\nabla (u^{\frac{2+\beta}{2}})|^2  + \frac{2}{1+\beta} \iint_{Q^{\tau}_R} \eta^2 u^{2+\beta} |\nabla \xi|^2
\end{split}
\end{equation*}
where we used the $\epsilon$-version of Cauchy inequality
\[
ab \le \epsilon a^2 + \frac{1}{4\epsilon} b^2, \quad \forall \epsilon>0
\]
in the second line by taking $\epsilon=\frac{1+\beta}{2}$.
Using the above three inequalities, we obtain
\begin{equation*}
\begin{split}
 &\frac{p}{p+1+\beta} \int_{B_R}  \xi^2   u^{p+1+\beta}(\tau) + \frac{2(1+\beta)}{(2+\beta)^2} \iint_{Q^{\tau}_R} \xi^2 \eta^2 |\nabla (u^{\frac{2+\beta}{2}})|^2 \\
 &\le \frac{C}{|R-r|^2}\iint_{Q^\tau_R}  u^{p+1+\beta} + \frac{2}{1+\beta} \iint_{Q^{\tau}_R} \eta^2 u^{2+\beta} |\nabla \xi|^2.
\end{split}
\end{equation*}
Taking the supremum over $\tau\in(-r^2,0)$ for the time-slice term and letting $\tau\uparrow0$ for the gradient term gives
\begin{equation}\label{eq:sup in t cut off R-r basic inequality}
\sup_{t\in(-r^2,0)}\int_{B_r}u^{p+1+\beta}(t)\,\ud x
+\iint_{Q_r}|\nabla u^{(2+\beta)/2}|^2
\le \frac{C(1+\beta)}{(R-r)^2}\iint_{Q_R}(u^{p+1+\beta}+u^{2+\beta}).
\end{equation}
Set $m_\beta=p+1+\beta$. For every $1/2<\rho\le2$, the inclusion $B_{1/2}\subset B_\rho$ and the hypothesis give
\[
\delta\le\|u(t)\|_{L^p(B_{1/2})}
\le\|u(t)\|_{L^p(B_\rho)}
\le |B_\rho|^{1/p-1/m_\beta}\|u(t)\|_{L^{m_\beta}(B_\rho)}.
\]
Consequently,
\be\label{eq:lower bound of L p+1+alpha}
\|u(t)\|_{L^{m_\beta}(B_\rho)}\ge c(n,p,\delta)>0,
\ee
uniformly in $\beta\ge0$ and $\rho\in(1/2,2]$. Since $m_\beta-(2+\beta)=p-1$, H\"older's inequality yields
\begin{equation}\label{eq:2+alpha bound by p+1+alpha}
\begin{split}
\int_{B_\rho}u^{2+\beta}
&\le |B_\rho|^{(p-1)/m_\beta}
       \|u\|_{L^{m_\beta}(B_\rho)}^{2+\beta}\\
&\le C(n,p,\delta)\int_{B_\rho}u^{m_\beta}.
\end{split}
\end{equation}
This estimate applies on both $B_R$ and $B_r$ in the following argument. Then it follows from \eqref{eq:sup in t cut off R-r basic inequality} that
\begin{equation}\label{eq:sup in t cut off R-r basic inequality without 2+alpha}
\begin{split}
 &  \sup_{ t \in (-r^2,0)}  \int_{B_r}   u^{p+1+\beta}(t) +  \iint_{Q_r}  |\nabla (u^{\frac{2+\beta}{2}})|^2  \le C\frac{1+\beta}{|R-r|^2}\iint_{Q_R}   u^{p+1+\beta}.
\end{split}
\end{equation}
For $\sigma>1$, employing \eqref{eq:2+alpha bound by p+1+alpha} and Sobolev inequality \eqref{eq:Sobolev inequality Br}, we have
\begin{equation*}
\begin{split}
\int_{B_r} u^{(2+\beta)\sigma} &= \int_{B_r} u^{2+\beta} u^{(2+\beta)(\sigma -1)}\\
& \le  \left(  \int_{B_r} u^{\frac{n}{n-2}(2+\beta)}       \right)^{\frac{n-2}{n}}  \left(   \int_{B_r}  u^{\frac{n}{2}(2+\beta)(\sigma -1) }    \right)^{\frac{2}{n}} \\
& \le C \left(    \int_{B_r} |\nabla (u^{\frac{2+\beta}{2}})|^2 +   \int_{B_r} u^{2+\beta}      \right)\left(   \int_{B_r}  u^{\frac{n}{2}(2+\beta)(\sigma -1) }    \right)^{\frac{2}{n}} \\
&\le C \left(    \int_{B_r} |\nabla (u^{\frac{2+\beta}{2}})|^2 +   \int_{B_r} u^{p+1+\beta}      \right)\left(   \int_{B_r}  u^{\frac{n}{2}(2+\beta)(\sigma -1) }    \right)^{\frac{2}{n}}.
\end{split}
\end{equation*}
Integrating the above inequality over the time interval $(-r^2,0)$ and using \eqref{eq:sup in t cut off R-r basic inequality without 2+alpha}, we obtain
\begin{equation}\label{eq: u^2+alpha sigma}
\begin{split}
\iint_{Q_r} u^{(2+\beta)\sigma} &  \le C \left(    \iint_{Q_r} |\nabla (u^{\frac{2+\beta}{2}})|^2 +   \iint_{Q_r} u^{p+1+\beta}      \right)\left( \sup_{t\in(-r^2,0)}  \int_{B_r}  u^{\frac{n}{2}(2+\beta)(\sigma -1) }  (t) \right)^{\frac{2}{n}} \\
&  \le C  \frac{1+\beta}{|R-r|^2}   \left(      \iint_{Q_R} u^{p+1+\beta}      \right)\left( \sup_{t\in(-r^2,0)}  \int_{B_r}  u^{\frac{n}{2}(2+\beta)(\sigma -1) }  (t) \right)^{\frac{2}{n}} .
\end{split}
\end{equation}
Let
\[
p+1+\beta = \frac{n}{2}(2+\beta)(\sigma -1),
\]
then
\[
\sigma=1+ \frac{2}{n} \cdot \frac{p+1+\beta}{2+\beta}.
\]
For $p$ satisfying \eqref{eq:sobolev subcritical exponent}, i.e.
\[
2<p+1<\frac{2n}{n-2}, \quad \mbox{if} \quad n \ge 3,
\]
we have that
\[
1+\frac{2}{n} \le \sigma \le 1+\frac{2}{n-2}, \quad \forall \  \beta \ge 0.
\]
Then it follows from \eqref{eq:sup in t cut off R-r basic inequality without 2+alpha} and \eqref{eq: u^2+alpha sigma} that
\begin{equation}\label{eq: u^2+alpha sigma bound by p+1+alpha}
\begin{split}
\iint_{Q_r} u^{(2+\beta)\sigma}   \le C \left( \frac{1+\beta}{|R-r|^2}  \right)^{1+\frac{2}{n}}  \left(      \iint_{Q_R} u^{p+1+\beta}      \right)^{1+\frac{2}{n}} .
\end{split}
\end{equation}
For $k=0,1,2,3,...$, we take ${\beta}_0=0$, and
\begin{equation*}
\begin{split}
\sigma_k &= 1+ \frac{2}{n} \cdot \frac{p+1+{\beta}_k}{2+{\beta}_k},\\
p+1+{\beta}_{k+1}&= (2+{\beta}_k)\sigma_k.
\end{split}
\end{equation*}
Put $q_k=p+1+\beta_k$, $a=1+2/n$, and $\ell_c=n(p-1)/2$. Then $q_{k+1}=a q_k-(p-1)$ and $q_k=\ell_c+a^k(p+1-\ell_c)$. For $1/2<r<R\le2$ set $R_k=r+2^{-k}(R-r)$ and $X_k=\iint_{Q_{R_k}}u^{q_k}$. Applying \eqref{eq: u^2+alpha sigma bound by p+1+alpha} gives
\[
X_{k+1}\le C\left(\frac{(1+\beta_k)4^{k+1}}{(R-r)^2}\right)^a X_k^a.
\]
Taking the $a^{-(k+1)}$ power and iterating gives a bounded product of the factors independent of $R-r$, since $\beta_k=O(a^k)$ and $\sum_{k\ge0}(k+1)a^{-k}<\infty$. Since $q_k/a^k\to p+1-\ell_c>0$, we obtain
\[
\|u\|_{L^\infty(Q_r)}\le \frac{C}{(R-r)^N}
\left(\iint_{Q_R}u^{p+1}\right)^{1/(p+1-\ell_c)},
\]
where $C$ depends on $n,p,\delta$ and $N=2a/((a-1)(p+1-\ell_c))$. The same argument applies to $v$. Taking $r=1$ and $R=2$ proves the theorem.
\epf

\begin{cor}\label{cor:L^infty local rescaled}
Let $n\ge3$, let $p$ satisfy \eqref{eq:sobolev subcritical exponent}, and let $1/2<r<R<2$. Suppose $(u,v)$ is a smooth nonnegative solution pair of the equations in Theorem~\ref{thm:L^infty local normalize} on $Q_R$, with $\min\{\|u(t)\|_{L^p(B_{1/2})},\|v(t)\|_{L^p(B_{1/2})}\}\ge\delta>0$ for $-R^2<t<0$. Then
\[
\|u\|_{L^\infty(Q_r)}+\|v\|_{L^\infty(Q_r)}
\le \frac{C}{(R-r)^N}\|(u,v)\|_{L^{p+1}(Q_R)}^\theta,
\]
with the constants from the preceding proof.
\end{cor}

\subsubsection{\texorpdfstring{Weak solutions for diagonal data}{Weak solutions for diagonal data}}

\begin{thm}\label{thm:short time existence of weak solutoin}
Let $p$ satisfy \eqref{eq:sobolev subcritical exponent} and let $u_0=v_0=w_0$, where $w_0\in H_0^1(\Omega)\cap L^\infty(\Omega)$ is nonnegative and not identically zero. Then \eqref{eq:fast diffusion system} admits a global nonnegative weak solution $(w,w)$ with zero Dirichlet data and initial trace $w^p(\cdot,0)=w_0^p$. In particular, $\|w(t)\|_{L^p(\Omega)}$ is bounded away from zero on a sufficiently short initial interval.
\end{thm}

\bpf
On the diagonal the equation is $\partial_t(w^p)=\Delta w-\alpha w$. Fix a finite time horizon and a uniform time step $h$, starting from the prescribed initial datum $w_0$. Given $w_k$, let $w_{k+1}$ minimize over $H_0^1(\Omega)$ the functional
\[
\int_\Omega \left[
 \frac{|z|^{p+1}}{(p+1)h}-\frac{(w_k)^p z}{h}
 +\frac{|\nabla z|^2+\alpha z^2}{2}
 \right]\,\ud x.
\]
Coercivity and weak lower semicontinuity give a minimizer. Replacing it by its absolute value cannot increase the functional. The Euler--Lagrange equation is
\[
\frac{(w_{k+1})^p-(w_k)^p}{h}
 =\Delta w_{k+1}-\alpha w_{k+1}.
\]
With $M=\|w_0\|_\infty$, testing with $(w_{k+1}-M)_+$ gives $0\le w_{k+1}\le M$ by induction. Put $A(z)=\int_\Omega (|\nabla z|^2+\alpha z^2)\,\ud x$ and $P_k=(w_k)^p$. Testing with $w_{k+1}-w_k$, and using the Lipschitz bound $L=pM^{p-1}$ for $s^p$ on $[0,M]$, gives
\[
\frac{\|P_{k+1}-P_k\|_2^2}{Lh}
 +\frac12\bigl[A(w_{k+1})-A(w_k)+A(w_{k+1}-w_k)\bigr]\le0.
\]
Summation bounds the $H_0^1$ norms of the iterates and the $L^2$ time derivative of the piecewise linear interpolant of $P_k$. The interpolants are also uniformly bounded in $H_0^1$ in space, since $s^p$ is Lipschitz on $[0,M]$. Compact spatial embedding and equicontinuity in time give convergence along a subsequence in $C([0,S];L^2(\Omega))$ on each finite horizon $S$. The piecewise constant powers have the same limit, because the largest $L^2$ jump is $O(h^{1/2})$. Taking the $p$th root gives strong space-time convergence of the piecewise constant iterates, whose gradients converge weakly along a subsequence. Passing to the limit in the discrete weak formulation gives the scalar equation with zero boundary trace. The convergence of the interpolants at $t=0$ gives $w^p(0)=w_0^p$. A diagonal subsequence over integer horizons gives the global solution. Continuity of $w^p$ in $L^1$ and $w_0\not\equiv0$ give the claimed short-time lower bound.
\epf

\brem
To continue a weak solution beyond a finite time, one needs an endpoint trace from which the equation can be solved again. The preceding theorem provides this for diagonal data. For general pairs, positivity before the endpoint is not enough to justify this step or to ensure simultaneous extinction. The convergence theorem therefore assumes common extinction. This regime is realized, for instance, by $u=v=(T-t)_+^{1/(p-1)}W$, where $(W,W)$ is the positive stationary pair constructed below.
\erem

\subsection{Stationary system} \label{subsec:ex stationary}
\begin{defn}\label{defn:weak solution steady}
We say that $(\Phi,\phi)$ is a pair of weak solutions of \eqref{eq:steady system} if
\be
\begin{split}
 &\Phi,\phi  \in H_0^1(\om),\\
  &\Phi,\phi \ge 0,
  \end{split}
\ee
satisfy
\be
\begin{split}
& \int_{\om}  \Big( \nabla \Phi \cdot \nabla \xi_1 + \alpha \phi       \xi_1 \Big) \, \ud x   = \frac{p}{p-1}\int_{\om}  \Phi^p \xi_1  \, \ud x,\\
 & \int_{\om}  \Big( \nabla \phi \cdot \nabla \xi_2 + \alpha \Phi       \xi_2 \Big) \, \ud x   = \frac{p}{p-1}\int_{\om}  \phi^p \xi_2  \, \ud x,
\end{split}
\ee
for any  $ \xi_1, \,  \xi_2 \in C^{\infty}_c(\om )$ .
\end{defn}

\brem \label{rem:bootstrap steady and positivity by maximum}
A standard  bootstrap technique shows that $(\Phi,\phi)$ is a pair of classical solutions if $p$ lies in the Sobolev subcritical  range \eqref{eq:sobolev subcritical exponent} (see \cite{Berger1977}, p. 50). If either component vanishes identically, its equation and $\alpha>0$ force the other component to vanish identically as well.
\erem

Next, we prove the existence of solutions of \eqref{eq:steady system}.
\begin{prop}\label{prop:existence steady system}
Let $p$ satisfy \eqref{eq:sobolev subcritical exponent}. Then \eqref{eq:steady system} admits a positive classical solution pair of the form $(W,W)$.
\end{prop}

\begin{proof}
For functions $f,g\in H_0^1(\om)$, $f,g \not\equiv 0$ on $ \om$, we define
\[
E_{p,\alpha}[f,g]:= \frac{ \int_{\om} |\nabla f|^2\,\ud x + \int_{\om} |\nabla g|^2\,\ud x  +\int_{\om} 2\alpha|fg|\,\ud x   }{  \left(\int_{ \om} |f|^{p+1}\,\ud x + \int_{ \om} |g|^{p+1}\,\ud x\right)^{\frac{2}{p+1}}}.
\]
\begin{equation}\label{eq:energyY}
Y^{\mathrm{diag}}_{p,\alpha}
=\inf_{0\ne w\in H_0^1(\Omega)} E_{p,\alpha}[w,w].
\end{equation}
On the diagonal, the quotient is
\[
E_{p,\alpha}[w,w]
=2^{1-2/(p+1)}
\frac{\int_\Omega (|\nabla w|^2+\alpha w^2)\,\ud x}
{\left(\int_\Omega |w|^{p+1}\,\ud x\right)^{2/(p+1)}}.
\]
Let $\mu=2^{2/(p+1)-1}Y^{\mathrm{diag}}_{p,\alpha}$. By the Sobolev inequality, $\mu>0$. Normalize a minimizing sequence by $\|w\|_{L^{p+1}}=1$. The numerator bounds the sequence in $H_0^1(\Omega)$, and the subcritical compact embedding preserves the normalization in the limit. Weak lower semicontinuity gives a minimizer $w$, which may be taken nonnegative by replacing it with $|w|$. Its Euler--Lagrange equation is
\[
-\Delta w+\alpha w=\mu w^p\quad\text{in }\Omega,\qquad w=0\quad\text{on }\partial\Omega.
\]
Elliptic bootstrapping gives classical regularity. Since $w\ge0$, $w\not\equiv0$ and $-\Delta w+\alpha w\ge0$, the scalar strong maximum principle gives $w>0$ in $\Omega$. Set $W=((p-1)\mu/p)^{1/(p-1)}w$. The pair $(W,W)$ satisfies both equations of \eqref{eq:steady system}.
\end{proof}

The diagonal construction gives one example: a positive separable solution with the same profile in both components. If a nonnegative, nontrivial stationary pair satisfies $\phi=k\Phi$ for a constant $k>0$, subtracting the first equation of \eqref{eq:steady system} from the second divided by $k$ gives
\[
\alpha(k^{-1}-k)\Phi
=\frac{p}{p-1}(k^{p-1}-1)\Phi^p.
\]
At any point where $\Phi>0$, the two sides have opposite signs for $k\ne1$. Thus proportional profiles must coincide.

\smallskip
\noindent\textit{Outlook.} Which profiles can occur, and how the evolution selects them, remain questions for further study.

\section{ Upper and lower integral bounds of the solutions of the diffusion system } \label{sec:integral}
The diagonal constructions give examples of the behavior we wish to study. We now consider general solution pairs under the hypotheses of Theorem~\ref{thm:integralbounds}. The coupled Rayleigh quotient lets us estimate the two components together.
\begin{proof}
By \eqref{eq:fast diffusion system}, we have
\begin{equation}\label{eq:dt to L^p+1}
\begin{split}
   &\frac{\ud}{\ud t} \left(  \int_{ \Omega} u^{p+1} \,\ud x+\int_{ \Omega} v^{p+1} \,\ud x \right)  \\
   =& -\frac{p+1}{p} \bigg(\int_{\Omega}|\nabla u|^{2}\,\ud x +\int_{\Omega}|\nabla v|^{2}\,\ud x+\int_{\Omega}2 \alpha u v \,\ud x \bigg)\le 0,
   \end{split}
\end{equation}
and
\be\label{eq:dt to dirichlet energy}
\begin{split}
&\frac{\ud}{\ud t} \left(\int_{\Omega}|\nabla u|^{2}\,\ud x +\int_{\Omega}|\nabla v|^{2}\,\ud x+\int_{\Omega}2 \alpha u v \,\ud x \right)\\
=& -\frac{2}{p} \left( \int_{\om}  \frac{\left(\Delta u -\alpha v\right)^2 }{u^{p-1}} \,\ud x +   \int_{\om}  \frac{\left(\Delta v -\alpha u\right)^2 }{v^{p-1}} \,\ud x             \right) \le 0.
\end{split}
\ee
Let
\[
I(t)= Q_\alpha[u(t),v(t)]= \frac{  \displaystyle \int_{\Omega}|\nabla u|^{2}\,\ud x    +\int_{\Omega}|\nabla v|^{2}\,\ud x   +\int_{\Omega}2\alpha uv\,\ud x       }{ \displaystyle\left(\int_{ \Omega}u^{p+1} \,\ud x  + \int_{ \Omega}v^{p+1} \,\ud x\right)^{\frac{2}{p+1}}   }  .
\]
Then we compute
\begin{align*}
\frac{\ud}{\ud t} I(t)&=\frac{2}{p} \left(\int_{ \Omega} u^{p+1} \,\ud x  +  \int_{ \Omega} v^{p+1} \,\ud x \right) ^{-\frac{2}{p+1}}\\
&\cdot \Bigg[\frac{\displaystyle\bigg(\int_{\Omega}|\nabla u|^{2}\,\ud x + \int_{\Omega}|\nabla v|^{2}\,\ud x +\int_{\Omega}2\alpha uv\,\ud x  \bigg)^2}{\displaystyle\int_{ \Omega} u^{p+1} \,\ud x  +  \int_{ \Omega} v^{p+1} \,\ud x }\\
&\qquad - \left( \int_{\om}  \frac{\left(\Delta u -\alpha v\right)^2 }{u^{p-1}} \,\ud x +  \int_{\om}  \frac{\left(\Delta v -\alpha u\right)^2 }{v^{p-1}} \,\ud x   \right)          \Bigg]\\
&=\frac{2}{p} \left(\int_{ \Omega} u^{p+1} \,\ud x  +  \int_{ \Omega} v^{p+1} \,\ud x \right) ^{-\frac{2}{p+1}} \cdot \left[    A +B    \right]
\end{align*}
where
\begin{align*}
A=&\frac{\displaystyle\bigg(\int_{\Omega}|\nabla u|^{2}\,\ud x + \int_{\Omega}|\nabla v|^{2}\,\ud x +\int_{\Omega}2\alpha uv\,\ud x  \bigg)^2}{\displaystyle\int_{ \Omega} u^{p+1} \,\ud x  +  \int_{ \Omega} v^{p+1} \,\ud x }\\
&-\frac{\displaystyle\bigg(\int_{\Omega}|\nabla u|^{2}\,\ud x+ \int_{\Omega}\alpha uv\,\ud x  \bigg)^2}{\displaystyle\int_{ \Omega} u^{p+1} \,\ud x  }-\frac{\displaystyle\bigg(\int_{\Omega}|\nabla v|^{2}\,\ud x+ \int_{\Omega}\alpha uv\,\ud x  \bigg)^2}{\displaystyle\int_{ \Omega} v^{p+1} \,\ud x  }
\end{align*}
and
\begin{align*}
B=&\frac{\displaystyle\bigg(\int_{\Omega}|\nabla u|^{2}\,\ud x +\int_{\Omega} \alpha uv\,\ud x  \bigg)^2}{\displaystyle\int_{ \Omega} u^{p+1} \,\ud x  } - \int_{ \Omega} \frac{\left(\Delta u -\alpha v\right)^2}{u^{p-1}}\,\ud x \\
&+ \frac{\displaystyle\bigg(\int_{\Omega}|\nabla v|^{2}\,\ud x+ \int_{\Omega} \alpha uv\,\ud x  \bigg)^2}{\displaystyle\int_{ \Omega} v^{p+1} \,\ud x  } -  \int_{ \Omega} \frac{\left(\Delta v -\alpha u\right)^2}{v^{p-1}}\,\ud x.
\end{align*}
Firstly, we denote
\[
A_1= \int_{\Omega}|\nabla u|^{2}\,\ud x +\int_{\Omega}\alpha uv\,\ud x   , \quad A_2= \int_{\Omega}|\nabla v|^{2}\,\ud x +\int_{\Omega}\alpha uv\,\ud x
\]
and
\[
B_1=\displaystyle\int_{ \Omega} u^{p+1} \,\ud x , \quad B_2= \displaystyle\int_{ \Omega} v^{p+1} \,\ud x.
\]
Then
\begin{align*}
A&= \frac{ \left( A_1+A_2 \right)^2}{B_1+B_2} - \frac{{A_1}^2}{B_1} - \frac{{A_2}^2}{B_2}\\
&=-\frac{1}{B_1+B_2} \cdot \left(   \frac{B_2}{B_1}         {A_1}^2    +   \frac{B_1}{B_2}         {A_2}^2  - 2  A_1 A_2 \right)\\
&=-\frac{1}{B_1+B_2} \cdot \left(   \sqrt{\frac{B_2}{B_1}          }{A_1}   -  \sqrt{\frac{B_1}{B_2} }         {A_2}   \right)^2 .
\end{align*}
Hence,
\[
A \le 0.
\]
Secondly, since $u=v=0$ on $\pa \Omega$, we have
\[  \int_{ \Omega} u (\Delta u -\alpha v) \,\ud x  =-\left( \int_{\Omega}|\nabla u|^{2}\,\ud x + \int_{\Omega} \alpha uv\,\ud x  \right)
\]
and
\[  \int_{ \Omega} v (\Delta v -\alpha u) \,\ud x  =-\left( \int_{\Omega}|\nabla v|^{2}\,\ud x + \int_{\Omega} \alpha uv\,\ud x  \right).
\]
 Applying the Cauchy-Schwarz inequality, we obtain:
\begin{equation*}
  \bigg(\int_{\Omega}|\nabla u|^{2}\,\ud x    +\int_{\Omega} \alpha uv \,\ud x     \bigg)^{2}  \leq  \int_{ \Omega}u^{p+1} \,\ud x    \int_{ \Omega} \frac{\left(\Delta u -\alpha v\right)^2}{u^{p-1}}\,\ud x
\end{equation*}
and
\begin{equation*}
  \bigg(\int_{\Omega}|\nabla v|^{2}\,\ud x    +\int_{\Omega} \alpha uv \,\ud x     \bigg)^{2}  \leq  \int_{ \Omega}v^{p+1} \,\ud x    \int_{ \Omega} \frac{\left(\Delta v -\alpha u\right)^2}{v^{p-1}}\,\ud x.
\end{equation*}
Hence,
\[
B \le 0.
\]
Therefore, by the nonpositivity of $A$ and $B$,
\begin{equation}\label{eq:monoticity of I(t)}
\frac{\ud}{\ud t} I(t)\le 0.
\end{equation}
This proves the monotonicity of Rayleigh quotient that $t \mapsto Q_\alpha[u(t),v(t)]$ is nonincreasing. The assumed initial convergence in $H_0^1$, together with Sobolev embedding, gives convergence of the quotient at $t=0$. Hence
\[
Q_\alpha[u(t),v(t)] \le Q_\alpha[u_0,v_0].
\]
 Together with \eqref{eq:dt to L^p+1}, we see
\[
\frac{p}{p+1} \frac{\ud}{\ud t}  \|  (u,v)   \|^{p+1}_{L^{p+1}(\om)} = -Q_\alpha[u(t),v(t)]  \|  (u,v)   \|^{2}_{L^{p+1}(\om)} \ge -Q_\alpha[u_0,v_0] \|  (u,v)   \|^{2}_{L^{p+1}(\om)}.
\]
So,
\[
\frac{p}{p-1} \frac{\ud}{\ud t}  \|  (u,v)   \|^{p-1}_{L^{p+1}(\om)}  \ge -Q_\alpha[u_0,v_0]
\]
We define
\begin{equation*}
  Z(t):=   \|  (u,v)   \|^{p-1}_{L^{p+1}(\om)}=  \bigg(  \displaystyle\int_{ \Omega} u^{p+1}(x,t) \,\ud x  + \displaystyle\int_{ \Omega} v^{p+1}(x,t) \,\ud x  \bigg)^{\frac{p-1}{p+1}},
\end{equation*}
then
\[
\frac{p}{p-1}Z(t) \ge \frac{p}{p-1}Z(0)- Q_\alpha[u_0,v_0]t.
\]
This yields
\[
T^* \ge \frac{p}{p-1} \frac{Z(0)}{Q_\alpha[u_0,v_0]} = \frac{p}{p-1}\frac{    \displaystyle  \int_{ \Omega}u_0^{p+1} \,\ud x  + \int_{ \Omega}v_0^{p+1} \,\ud x     }{ \displaystyle \int_{\Omega}|\nabla u_0|^{2}\,\ud x    +\int_{\Omega}|\nabla v_0|^{2}\,\ud x   +\int_{\Omega}2 \alpha u_0v_0\,\ud x   },
\]
by letting $t \to T^*$. Next, it follows from \eqref{eq:dt to L^p+1} that
\begin{equation}\label{eq:Zt}
  Z'(t)= -\frac{p-1}{p}  I(t).
\end{equation}
In particular, if $p$ lies in the subcritical range \eqref{eq:sobolev subcritical exponent}, by Sobolev inequality \eqref{eq:sobolev-poincare} we find that for any nontrivial $u,v \in H_0^1(\om), u,v \ge 0$,
\be \label{eq:lower bound of It}
I(t) \ge 1/C
\ee
Such positive constant $C$ depends only on $n,p$ and $\Omega$. Subsequently, it follows from \eqref{eq:lower bound of It} and \eqref{eq:Zt} that
\[
  Z'(t)\le -C
\]
for some $C>0$ depending only on $n,p$ and $\Omega$. Then
\[
Z(T)-Z(t)\le -C(T-t).
\]
Taking $\liminf_{T\to (T^*)^-}$ on both sides and using the definition of $T^*$ \eqref{eq:extinction surmise}, we obtain
\begin{equation*}
Z(t)\ge C(T^*-t).
\end{equation*}
Hence, the first inequality follows.

By \eqref{eq:monoticity of I(t)} and \eqref{eq:Zt}, $Z'(t)$ is nondecreasing, so $Z$ is convex. Hence, for $0<s<t<T<T^*$, we have
\[
Z(t)\le Z(T)+\frac{Z(s)-Z(T)}{s-T} (t-T).
\]
Taking $\liminf_{T\to (T^*)^-}$ and $\lim_{s\to 0^+}$ on both sides, and using \eqref{eq:extinction surmise}, we have
\[
Z(t)\le \frac{Z(0)}{T^*} (T^*-t).
\]
Hence, the second inequality follows.

\end{proof}

\brem \label{rem:tilde u,v L p+1}
Let
\begin{equation} \label{eq:changeofvariable tilde}
 \left\{
\begin{array}{lll}
\tilde{u}(x,t)&=u(x,t)/(T^*-t)^{\frac{1}{p-1}} \ , \\
\tilde{v}(x,t)&=v(x,t)/(T^*-t)^{\frac{1}{p-1}} \ ,   \quad
\end{array}
\right.
\end{equation}
Equivalently,
\begin{equation*}
 \left\{
\begin{array}{lll}
\tilde{u}(x,t)&=U(x,\tau) \ \mbox{ with }\  t=T^*(1-e^{-\tau}), \\
\tilde{v}(x,t)&=V(x,\tau) \ \mbox{ with }\  t=T^*(1-e^{-\tau}).  \quad
\end{array}
\right.
\end{equation*}
Let $(u,v)$ satisfy the hypotheses of Theorem~\ref{thm:integralbounds}, including the common extinction hypothesis, and define $(\tilde u,\tilde v)$ as above. Then
\[
\int_{ \Omega} \tilde{u}^{p+1} \,\ud x+\int_{ \Omega} \tilde{v}^{p+1} \,\ud x
\]
is nonincreasing in $t$. Thus,
\[
\int_{ \Omega} U^{p+1} \,\ud x+\int_{ \Omega} V^{p+1} \,\ud x
\]
is nonincreasing in $\tau$.
\erem
\bpf
We define
\begin{equation*}
  \tilde{Z}(t):=   \|  (\tilde{u},\tilde{v})   \|^{p-1}_{L^{p+1}(\om)}=  \bigg(  \displaystyle\int_{ \Omega} \tilde{u}^{p+1}(x,t) \,\ud x  + \displaystyle\int_{ \Omega} \tilde{v}^{p+1}(x,t) \,\ud x  \bigg)^{\frac{p-1}{p+1}},
\end{equation*}
then
\[
\tilde{Z}(t) = \frac{Z(t)}{T^*-t}.
\]
From the proof of Theorem \ref{thm:integralbounds}, we know that $Z(t)$ is convex, and common extinction gives $Z(T^*)=0$ by continuous extension. Therefore,
\[
 -\tilde{Z}(t)= \frac{Z(T^*)-Z(t)}{T^*-t}
\]
is nondecreasing for $t \in (0,T^*)$.
\epf

\section{Stability of the separable solutions: convergence in a sequence of times}\label{sec:con}

In this section, we prove Theorem \ref{thm:convergence}.

\bpf
Define the functional
\[
F(f,g):= \int_{\om} \left[ \frac{1}{2}|\nabla f|^2    + \frac{1}{2}|\nabla g|^2  +\alpha fg -\frac{p}{(p-1)(p+1)}    \left(  f^{p+1} + g^{p+1} \right) \right] \, \ud x.
\]
Let
\[
g(\tau)=F(U(\cdot, \tau),V(\cdot, \tau)),
\]
where $(U,V)$ satisfies \eqref{eq:fast diffusion system scale}. Differentiating gives
\be \label{eq:g't}
\begin{split}
g'(\tau) &= \left\langle  \pa_U F(U,V) , \pa_\tau U \right\rangle +   \left\langle  \pa_V F(U,V) , \pa_\tau V \right\rangle  \\
&=\int_\om \left[ \left(  -\Delta U +\alpha V -\frac{p}{p-1} U^p        \right) \pa_\tau U \right] \, \ud x + \int_\om \left[ \left(  -\Delta V +\alpha U -\frac{p}{p-1} V^p        \right) \pa_\tau V \right] \, \ud x \\
&= - p \left(     \int_\om  U^{p-1} |\pa_\tau U|^2     \, \ud x     +       \int_\om  V^{p-1} |\pa_\tau V|^2     \, \ud x       \right)\\
& \le 0,
\end{split}
\ee
since $U,V$ are both nonnegative. From Theorem \ref{thm:integralbounds}, we have \eqref{eq:integral bounds rescaled}. It shows that
\[
 \int_{\om} U^{p+1}(x,\tau)\,\ud x+           \int_{\om} V^{p+1}(x,\tau)\,\ud x
\]
is bounded in $\tau$, which means $g(\tau)$ is bounded below. Therefore
\[
\lim_{\tau \to \infty} g(\tau)
\]
exists and there is a sequence of times
\[
\tau_n \to \infty   \quad  \mbox{as} \quad n \to \infty
\]
 such that
 \[
\lim_{n \to \infty} g'(\tau_n) = 0.
 \]
We shall show for a subsequence of $\{ \tau_n \}_{n=1}^\infty$, again labeled as $\tau_n$, that
\begin{equation*}
\begin{split}
U(\cdot, \tau_n) \to \Phi(\cdot)  \\
  V(\cdot, \tau_n) \to \phi(\cdot)
\end{split}
\end{equation*}
as $n \to \infty$ in $\mathring{W}^{1,\frac{p+1}{p}}(\om)$. To do this, notice that
\be \label{eq: Delta U estimate}
\begin{split}
&\int_{\om} |\Delta U(x,\tau_n)|^{\frac{p+1}{p}} \, \ud x\\
&\le C(p,\alpha) \left(     \int_{\om} |U^p(x,\tau_n)|^{\frac{p+1}{p}} \, \ud x   +  \int_{\om} |V(x,\tau_n)|^{\frac{p+1}{p}} \, \ud x         +    \int_{\om}  |U^{p-1}(x,\tau_n)\pa_\tau U(x,\tau_n) |^{\frac{p+1}{p}}  \, \ud x   \right)
\end{split}
\ee
and
\be \label{eq: Delta V estimate}
\begin{split}
&\int_{\om} |\Delta V(x,\tau_n)|^{\frac{p+1}{p}} \, \ud x\\
&\le C(p,\alpha) \left(     \int_{\om} |V^p(x,\tau_n)|^{\frac{p+1}{p}} \, \ud x   +  \int_{\om} |U(x,\tau_n)|^{\frac{p+1}{p}} \, \ud x         +    \int_{\om}  |V^{p-1}(x,\tau_n)\pa_\tau V(x,\tau_n) |^{\frac{p+1}{p}}  \, \ud x   \right),
\end{split}
\ee
since $(U,V)$ satisfies \eqref{eq:fast diffusion system scale}. For the last terms in \eqref{eq: Delta U estimate} and \eqref{eq: Delta V estimate}, we have
\be \label{eq: Delta U estimate-2}
\begin{split}
 & \int_{\om}  |U^{p-1}(x,\tau_n)\pa_\tau U(x,\tau_n) |^{\frac{p+1}{p}}  \, \ud x   \\
 &\le  \left(  \int_{\om}  |U(x,\tau_n)|^{p+1}  \, \ud x  \right)^{\frac{p-1}{2p}}\left( \int_{\om}  |U(x,\tau_n) |^{p-1}   |\pa_\tau U(x,\tau_n) |^{2}  \, \ud x    \right)^{\frac{p+1}{2p}}\\
 &=  \left(  \int_{\om}  U^{p+1}(x,\tau_n)  \, \ud x  \right)^{\frac{p-1}{2p}}\left( \int_{\om}  U^{p-1}(x,\tau_n)    |\pa_\tau U(x,\tau_n) |^{2}  \, \ud x    \right)^{\frac{p+1}{2p}}\\
 &\le C \|U(\cdot,\tau_n)  \|^{\frac{(p+1)(p-1)}{2p}}_{L^{p+1}(\om)} |g'(\tau_n)|^{\frac{p+1}{2p}}.
\end{split}
\ee
Similarly,
\be \label{eq: Delta V estimate-2}
\begin{split}
 & \int_{\om}  |V^{p-1}(x,\tau_n)\pa_\tau V(x,\tau_n) |^{\frac{p+1}{p}}  \, \ud x  \le C \|V(\cdot,\tau_n)  \|^{\frac{(p+1)(p-1)}{2p}}_{L^{p+1}(\om)} |g'(\tau_n)|^{\frac{p+1}{2p}}.
\end{split}
\ee
Combining \eqref{eq: Delta U estimate}, \eqref{eq: Delta V estimate}, \eqref{eq: Delta U estimate-2} and \eqref{eq: Delta V estimate-2}, we obtain
\be \label{eq: Delta U and Delta V estimate}
\begin{split}
&\int_{\om} |\Delta U(x,\tau_n)|^{\frac{p+1}{p}} \, \ud x +\int_{\om} |\Delta V(x,\tau_n)|^{\frac{p+1}{p}} \, \ud x \\
&\le C  \| \left( U(\cdot,\tau_n),V(\cdot,\tau_n) \right)  \|^{\frac{(p+1)(p-1)}{2p}}_{L^{p+1}(\om)} |g'(\tau_n)|^{\frac{p+1}{2p}}  + C\| \left( U(\cdot,\tau_n),V(\cdot,\tau_n) \right)  \|^{p+1}_{L^{p+1}(\om)} \\
& \ \ + C(p,\alpha)|\om|^{\frac{p-1}{p}}\| \left( U(\cdot,\tau_n),V(\cdot,\tau_n) \right)  \|_{L^{p+1}(\om)}^{\frac{p+1}{p}},
\end{split}
\ee
where we used
\[
\int_{\om} |U(x,\tau_n)|^{\frac{p+1}{p}} \, \ud x \le  |\om|^{\frac{p-1}{p}}  \left(  \int_{\om} |U(x,\tau_n)|^{p+1} \, \ud x \right)^{\frac{1}{p}}
\]
and
\[
\int_{\om} |V(x,\tau_n)|^{\frac{p+1}{p}} \, \ud x \le  |\om|^{\frac{p-1}{p}}  \left(  \int_{\om} |V(x,\tau_n)|^{p+1} \, \ud x \right)^{\frac{1}{p}}.
\]
By $g'(\tau_n) \to 0$ and \eqref{eq:integral bounds rescaled},
\[
\int_{\om} |\Delta U(x,\tau_n)|^{\frac{p+1}{p}} \, \ud x +\int_{\om} |\Delta V(x,\tau_n)|^{\frac{p+1}{p}} \, \ud x
\]
is bounded in $\{ \tau_n \}_{n=1}^\infty$. Put $q=(p+1)/p\in(1,2)$. Since $\Omega$ is smooth and both functions have zero Dirichlet boundary values, the global elliptic estimate gives
\[
\begin{split}
\|U(\cdot,\tau_n)\|_{W^{2,q}(\Omega)}
+\|V(\cdot,\tau_n)\|_{W^{2,q}(\Omega)}
\le C_{\Omega,q}\big(&\|\Delta U(\cdot,\tau_n)\|_{L^q(\Omega)}
+\|\Delta V(\cdot,\tau_n)\|_{L^q(\Omega)}\\
&+\|U(\cdot,\tau_n)\|_{L^q(\Omega)}
+\|V(\cdot,\tau_n)\|_{L^q(\Omega)}\big).
\end{split}
\]
The Laplacian terms have just been bounded. The remaining terms are controlled by \eqref{eq:integral bounds rescaled}, since $q<p+1$ and $\Omega$ is bounded. Thus the sequence is bounded in $W^{2,q}(\Omega)\cap W_0^{1,q}(\Omega)$, and the compact embedding into $W^{1,q}(\Omega)$ applies.
Hence there exists a subsequence, again labeled by $U(\cdot, \tau_n), V(\cdot, \tau_n)$ such that
\begin{equation} \label{eq:U V tend in W space}
\begin{split}
U(\cdot, \tau_n) \to \Phi(\cdot)   \\
  V(\cdot, \tau_n) \to \phi(\cdot)
\end{split}
\end{equation}
as $n \to \infty$ in $\mathring{W}^{1,\frac{p+1}{p}}(\om)$ for some functions $\Phi(\cdot)$ and $\phi(\cdot)$.
Moreover, since $g(\tau)$ and $\|(U,V)\|_{L^{p+1}(\om)}$ are bounded, we see 
\[
\int_\om |\nabla U(x,\tau)|^2 \, \ud x +\int_\om |\nabla V(x,\tau)|^2 \, \ud x
\]
is bounded, which gives that $ \|U(x,\tau)\|_{H_0^1(\om)} +\|V(x,\tau)\|_{H_0^1(\om)}  $ is uniformly bounded in $\tau$. It implies that $\Phi(\cdot)$ and $\phi(\cdot)$ must be in
$H_0^1(\om)$ and
\begin{equation}\label{eq:weak convergence in H}
\begin{split}
U(\cdot, \tau_n) \rightharpoonup \Phi(\cdot)   \\
  V(\cdot, \tau_n) \rightharpoonup \phi(\cdot)
\end{split}
\end{equation}
as $n\to\infty$ weakly in $H_0^1(\om)$. The subcritical compact embedding gives strong convergence in $L^{p+1}(\Omega)$ after passing to a further subsequence. In particular,
\begin{equation}\label{eq: UV integral p+1}
\begin{split}
\int_{\om} U^{p+1}(\cdot, \tau_n) \, \ud x  \to\int_{\om}  \Phi^{p+1}(\cdot) \, \ud x   \\
\int_{\om} V^{p+1}(\cdot, \tau_n) \, \ud x  \to\int_{\om}  \phi^{p+1}(\cdot) \, \ud x
\end{split}
\end{equation}
as $n \to \infty$, since $2 < p+1 < \frac{2n}{n-2}$ if $n > 2$ and $p+1 > 2$ if $n < 3$. Also, by Remark \ref{rem:tilde u,v L p+1}, we know that
$\int_\om U^{p+1} \, \ud x + \int_\om V^{p+1} \, \ud x$ is nonincreasing; hence
\begin{equation*}
\begin{split}
\int_{\om} U^{p+1}(\cdot, \tau) \, \ud x +  \int_{\om} V^{p+1}(\cdot, \tau) \, \ud x  \to    \int_{\om}  \Phi^{p+1}(\cdot) \, \ud x  + \int_{\om}  \phi^{p+1}(\cdot) \, \ud x
\end{split}
\end{equation*}
as $\tau \to \infty$.
Next, we show that $(\Phi(\cdot),\phi(\cdot))$ is a pair of weak solutions to \eqref{eq:steady system}. To do this, we multiply equations in \eqref{eq:fast diffusion system scale} by any $\xi_1(x),\xi_2(x) \in C^{\infty}_c(\om)$ respectively and integrate them over $\om$. Then integration by parts gives
\begin{equation}\label{eq: fast diffusion system scale test xi1 xi2}
\begin{split}
p\int_{\om} \xi_1 U^{p-1}\pa_\tau U \, \ud x = -\int_{\om} \left(\nabla U \cdot \nabla \xi_1 + \alpha V \xi_1 - \frac{p}{p-1}U^p \xi_1                      \right)\, \ud x, \\
p\int_{\om} \xi_2 V^{p-1}\pa_\tau V \, \ud x = -\int_{\om} \left(\nabla V \cdot \nabla \xi_2 + \alpha U \xi_2 - \frac{p}{p-1}V^p \xi_2                      \right)\, \ud x.
\end{split}
\end{equation}
The left hand sides of \eqref{eq: fast diffusion system scale test xi1 xi2} tend to zero as $\tau_n \to \infty$, since
\begin{equation*}
\begin{split}
&p\left|\int_{\om} \xi_1 U^{p-1}\pa_\tau U \, \ud x \right| + p\left|\int_{\om} \xi_2 V^{p-1}\pa_\tau V \, \ud x\right| \\
&\le C \left(  \max_\om |\xi_1|  \left( \int_{\om} U^{p-1}(x,\tau_n)\,\ud x \right)^{\frac{1}{2}}  + \max_\om |\xi_2|  \left( \int_{\om} V^{p-1}(x,\tau_n)\,\ud x \right)^{\frac{1}{2}}\right)|g'(\tau_n)|^{\frac{1}{2}} \to 0
\end{split}
\end{equation*}
as $n \to \infty$. Therefore, with $\tau_n \to \infty$, \eqref{eq: fast diffusion system scale test xi1 xi2} become
\begin{equation}\label{eq: fast diffusion system scale test xi1 xi2 becomes steady}
\begin{split}
\int_{\om} \left(\nabla \Phi \cdot \nabla \xi_1 + \alpha \phi \xi_1 - \frac{p}{p-1}\Phi^p \xi_1                      \right)\, \ud x = 0, \\
\int_{\om} \left(\nabla \phi \cdot \nabla \xi_2 + \alpha \Phi \xi_2 - \frac{p}{p-1}\phi^p \xi_2                      \right)\, \ud x=0,
\end{split}
\end{equation}
by weak $H_0^1$ convergence and strong $L^{p+1}$ convergence. Indeed, $U(\tau_n)^p\to\Phi^p$ and $V(\tau_n)^p\to\phi^p$ in $L^{(p+1)/p}(\Omega)$. Thus \eqref{eq: fast diffusion system scale test xi1 xi2 becomes steady} shows that $(\Phi,\phi)$ is a weak solution pair of \eqref{eq:steady system}.

The limits $\Phi$ and $\phi$ are nonnegative by strong convergence in $L^{p+1}(\Omega)$. The lower bound in \eqref{eq:integral bounds rescaled} passes to the limit, so the pair is nontrivial. By Remark \ref{rem:bootstrap steady and positivity by maximum}, both components are individually nontrivial, and elliptic bootstrapping in the Sobolev subcritical range gives classical regularity. Strong convergence in $L^{p+1}(\Omega)$ also gives
\be \label{eq:UV integral limit}
\lim_{n\to \infty} \int_{\om} U(x,\tau_n)V(x,\tau_n) \, \ud x = \int_{\om} \Phi(x)\phi(x)     \, \ud x.
\ee

Finally, we establish that
\begin{equation} \label{eq:U V strong convergence in H space}
\begin{split}
U(\cdot, \tau_n) \to \Phi(\cdot)   \\
  V(\cdot, \tau_n) \to \phi(\cdot)
\end{split}
\end{equation}
as $n \to \infty$ in $H_0^1(\om)$. To do this, we multiply equations of \eqref{eq:fast diffusion system scale} by $U$ and $V$ respectively, and integrate by parts to get
\begin{equation}\label{eq: fast diffusion system scale test U V}
\begin{split}
p\int_{\om}  U^{p}\pa_\tau U \, \ud x = -\int_{\om} \left( |\nabla U|^2  + \alpha UV  - \frac{p}{p-1}U^{p+1}                      \right)\, \ud x, \\
p\int_{\om}  V^{p}\pa_\tau V \, \ud x = -\int_{\om} \left(|\nabla V|^2  + \alpha U V - \frac{p}{p-1}V^{p+1}                       \right)\, \ud x.
\end{split}
\end{equation}
The weighted Cauchy--Schwarz inequality and \eqref{eq:g't} give
\[
\left|\int_\Omega U^p U_\tau\right|
\le \left(\int_\Omega U^{p+1}\right)^{1/2}
     \left(\int_\Omega U^{p-1}|U_\tau|^2\right)^{1/2}
\longrightarrow0
\]
along $\tau_n$, and the same holds for $V$. Thus
\begin{equation*}
\begin{split}
\lim_{n\to \infty} \int_{\om} \bigg[ |\nabla U(x,\tau_n)|^2 +  \alpha U(x,\tau_n)V(x,\tau_n) \bigg]  \, \ud x = \lim_{n\to \infty}\int_{\om} \frac{p}{p-1}U^{p+1}(x,\tau_n)   \, \ud x,\\
\lim_{n\to \infty} \int_{\om} \bigg[ |\nabla V(x,\tau_n)|^2 +  \alpha U(x,\tau_n)V(x,\tau_n) \bigg]  \, \ud x = \lim_{n\to \infty}\int_{\om} \frac{p}{p-1}V^{p+1}(x,\tau_n)   \, \ud x.
\end{split}
\end{equation*}
From \eqref{eq: UV integral p+1} and \eqref{eq:UV integral limit},
\begin{equation*}
\begin{split}
\lim_{n\to \infty} \int_{\om}  |\nabla U(x,\tau_n)|^2  \, \ud x + \int_{\om}  \alpha \Phi(x)\phi(x)  \, \ud x   = \int_{\om} \frac{p}{p-1}\Phi^{p+1}(x)   \, \ud x,\\
\lim_{n\to \infty} \int_{\om}  |\nabla V(x,\tau_n)|^2  \, \ud x + \int_{\om}  \alpha \Phi(x)\phi(x)  \, \ud x   = \int_{\om} \frac{p}{p-1}\phi^{p+1}(x)   \, \ud x.
\end{split}
\end{equation*}
Therefore,
\begin{equation} \label{eq:U V H space norm tend }
\begin{split}
\| U(\cdot, \tau_n) \|_{H_0^1(\om)} \to \| \Phi(\cdot) \|_{H_0^1(\om)}  \\
 \| V(\cdot, \tau_n) \|_{H_0^1(\om)} \to\| \phi(\cdot) \|_{H_0^1(\om)}
\end{split}
\end{equation}
as $n \to \infty$, since $(\Phi,\phi)$ satisfies \eqref{eq:steady system}, and hence
\begin{equation*}
\begin{split}
 \int_{\om}  |\nabla \Phi(x)|^2  \, \ud x + \int_{\om}  \alpha \Phi(x)\phi(x)  \, \ud x   = \int_{\om} \frac{p}{p-1}\Phi^{p+1}(x)   \, \ud x,\\
 \int_{\om}  |\nabla \phi(x)|^2  \, \ud x + \int_{\om}  \alpha \Phi(x)\phi(x)  \, \ud x   = \int_{\om} \frac{p}{p-1}\phi^{p+1}(x)   \, \ud x.
\end{split}
\end{equation*}
Since the norms converge and we have already established the weak convergence \eqref{eq:weak convergence in H}, we then in fact have the strong convergence \eqref{eq:U V strong convergence in H space}.

\epf

\bigskip

\noindent\textbf{Declarations of interest}: None.

\bigskip

\noindent \textbf{Acknowledgement:}
The author would like to thank his father for his valuable support and encouragement.

\end{document}